\documentclass[a4paper]{article}
\usepackage{authblk}
\usepackage{amsmath}
\usepackage{amsfonts}
\usepackage{amssymb}
\usepackage{amsthm}
\usepackage{bbm}
\usepackage{graphicx}
\usepackage[utf8]{inputenc}
\usepackage[square,numbers,sort]{natbib}
\usepackage{colonequals}
\usepackage{overpic}
\usepackage{booktabs}
\usepackage{todonotes}
\usepackage{enumitem}
\usepackage{tabularx}
\usepackage{hyperref}

\DeclareMathOperator{\tr}{tr}
\newcommand{\R}{\mathbb{R}}

\DeclareMathOperator{\dist}{dist}

\DeclareMathOperator{\sym}{sym}
\DeclareMathOperator{\SOdrei}{SO(3)}

\newcommand{\parder}[2]{\frac{\partial #1}{\partial #2}}
\newcommand{\totalder}[2]{\frac{\text{d} #1}{\text{d} #2}}
\providecommand{\abs}[1]{{\lvert#1\rvert}}
\providecommand{\norm}[1]{\lVert#1\rVert}
\newcommand{\identity}{\text{Id}}
\newcommand{\inj}{\operatorname{inj}}
\newcommand{\hessian}{\operatorname{Hess}}

\DeclareMathOperator*{\argmin}{arg\,min}

\newtheorem{theorem}{Theorem}
\newtheorem{lemma}[theorem]{Lemma}
\newtheorem{definition}{Definition}

\theoremstyle{remark}
\newtheorem{remark}{Remark}
\newtheorem{problem}{Problem}
\newtheorem*{example}{\rm \textbf{Example}}

\graphicspath{{gfx/}}

\title{Numerical Treatment of a Geometrically Nonlinear Planar Cosserat Shell Model}
\author{Oliver Sander%
        \footnote{Corresponding author: Oliver Sander, Institut für Geometrie und Praktische Mathematik, RWTH Aachen,
                  Templergraben 55, 52056 Aachen, Germany, email: \href{mailto:sander@igpm.rwth-aachen.de}{sander@igpm.rwth-aachen.de}
                 },\;
        Patrizio Neff%
        \footnote{Patrizio Neff, Head of Lehrstuhl für Nichtlineare Analysis und Modellierung, Fakultät für Mathematik,
                  Universität Duisburg--Essen, Campus Essen, Thea-Leymann Str. 9, 45127 Essen, Germany, email:
                  \href{mailto:patrizio.neff@uni-due.de}{patrizio.neff@uni-due.de}, Tel.:~+49-201-183-4243
                 }
        \; and Mircea B\^{\i}rsan%
        \footnote{Mircea B\^{\i}rsan, Lehrstuhl für Nichtlineare Analysis und Modellierung, Fakultät für Mathematik,
                  Universität Duisburg--Essen, Campus Essen, Thea-Leymann Str. 9, 45127 Essen, Germany, email:
                  \href{mailto:mircea.birsan@uni-due.de}{mircea.birsan@uni-due.de}; and Department of Mathematics, University ``A.I. Cuza'' of Ia\c si, 700506
                  Ia\c si, Romania
                 }
       }

\begin{document}

\maketitle

\begin{abstract}
We present a new way to discretize a geometrically nonlinear elastic planar Cosserat shell. The kinematical model is similar  to the general 6-parameter resultant shell model with drilling rotations.  The discretization uses geodesic finite elements, which leads
to an objective discrete model which naturally allows arbitrarily large rotations.  Finite elements of any
approximation order can be constructed.  The resulting algebraic problem is a minimization problem posed
on a nonlinear finite-dimensional Riemannian manifold.  We solve this problem using a Riemannian
trust-region method, which is a generalization of Newton's method that converges globally
without intermediate loading steps.  We present the continuous model and the discretization,
discuss the properties of the discrete model, and show several numerical examples, including wrinkles of thin elastic sheets in shear.
\end{abstract}

\section{Introduction}

We consider the numerical treatment of a geometrically nonlinear hyperelastic planar Cosserat shell model.  This model has been obtained
by dimensional reduction from a full three-dimensional Cosserat continuum model.  Its degrees of freedom are
the displacement $m$ of the shell midsurface, together with the orientation of an orthonormal director triple
$\overline{R}$ at each point.  Consequently, if $\omega$ denotes the two-dimensional parameter domain, configurations
of such a shell are pairs of functions
\begin{equation*}
 (m,\overline{R}) : \omega \to \R^3 \times \text{SO(3)},
\end{equation*}
of suitable smoothness.
We consider a hyperelastic material law of the form
\begin{equation}
\label{eq:general_energy}
  I(m,\overline{R})
  =
  \int_\omega h\,W_{\text{mp}}(\overline{U})+ h\,W_{\text{curv}}(\mathfrak{K}_s)+ \frac{h^3}{12}\, W_{\text{bend}}(\mathfrak{K}_b)\,\mathrm{d}\omega + \text{external loads},
\end{equation}
where $W_\text{mp}$ is the membrane energy, $W_\text{bend}$ is the bending energy, and $W_\text{curv}$ is a curvature
term depending only on the orientation field $\overline{R}$.  This energy, originally proposed in \cite{Neff_plate04_cmt,Neff_plate07_m3as},
is a second-order model, frame-invariant, and allows for large elastic strains and finite rotations.  The membrane
contribution $W_\text{mp}$ is polyconvex, and uniformly Legendre--Hadamard-elliptic.
Existence of minimizers in the space $H^1(\omega,\R^3) \times W^{1,q}(\omega,\text{SO(3)})$ has been shown in~\cite{Neff_plate04_cmt,Neff_plate07_m3as} for any $q\geq 2$.

In this article we consider planar shells only, i.e., we assume that the undeformed configuration $(m_0, \overline{R}_0) : (x,y) \mapsto ((x,y,0),\text{Id})$
is a stress-free state. However, our numerical treatment can also be generalized to a general nonplanar shell model.
We arrive at the planar model in two steps:  First, dimensional reduction of a parent three-dimensional Cosserat model yields a
shell model with a quadratic membrane energy, suitable only for small membrane strains.  We then generalize
this shell model to obtain the finite-strain membrane term.

The shell formulation presented here is closely related to the theory of 6-parameter shells with drilling rotations
\cite{Pietraszkiewicz-book04,Libai98,Eremeyev06}. A detailed comparison between the two approaches in the case of plates has been given in
\cite{Birsan-Neff-AnnRom-2012,Birsan-Neff-JElast-2013}, where we have shown  existence results for isotropic, orthotropic, and composite plates. In~\cite{Birsan-Neff-MMS-2014} we have adapted the  methods of \cite{Neff_plate04_cmt} to
prove the existence of minimizers for geometrically nonlinear 6-parameter shells. In~\cite{Birsan-Neff-L54-2014} we have considered shells insensitive to drilling rotations, and established a useful representation theorem for this case (corresponding to the Cosserat couple modulus $\mu_c=0$).\smallskip

We also mention that the kinematic assumption underlying our Cosserat shell formulation is similar to
the one used in describing a viscoelastic membrane or a viscoelastic rod, see \cite{Neff_membrane_plate03,Neff_membrane_existence03,Neff_Beyrouthy_11,Neff_membrane_Weinberg07}.
Indeed, the viscoelastic membrane is based on the same kinematics, but the independent rotations are
evolving through a local evolution equation, whereas for the Cosserat planar shell  model, they
are determined by energy minimization.

\bigskip

Problems with directional or orientational degrees of freedom are notoriously difficult to discretize.  This difficulty
is caused by the nonlinearity of the orientation configuration space $W^{1,q}(\omega,\text{SO(3)})$ (or, in fact,
any space of functions mapping into $\text{SO(3)}$).  As a consequence, discretization methods based on piecewise linear
or piecewise polynomial functions cannot be formulated directly for such spaces.  Instead, previous discretizations have used
ad hoc approaches, each with its particular shortcomings.

An obvious approach uses Euler angles to describe the rotations, and finite elements to discretize
the angles~\cite{wriggers_gruttmann:1993}.
However, this leads to instabilities near certain configurations, and such models are
suitable only for situations with moderately large rotations~\cite{gruttmann_wagner_meyer_wriggers:1993}.
Also, the resulting discrete models are generally not objective.

Alternatively, rotations can be interpolated by means of the Lie algebra $\mathfrak{so}(3)$, i.e., the tangent space at the identity rotation.
A rotation $R \in \text{SO}(3)$ is represented as a rotation vector $a \in \mathfrak{so}(3)$
with $R = \exp a$.  Since $\mathfrak{so}(3)$ is a linear space, the rotation
vectors $a$ can be interpolated normally using finite elements of first or higher order%
~\cite{mueller:2009,muench:2007}.
This approach works only for orientation values bounded away from the cut locus
of the identity rotation.  To deal with larger rotations, \cite{muench:2007} switches to a
different tangent space when large rotations are detected.

Unfortunately, using a fixed tangent space for interpolation introduces a preferred direction
into the discrete model.  The discrete solution therefore depends on the orientation of the observer,
and objectivity is not preserved.

For their model of a shell with a single director, \citet{simo_fox:1989} propose to avoid nonlinear
interpolation altogether.  Instead, they introduce the director vector directions at the quadrature points
as separate variables~\cite{simo_fox_rifai:1990}.  The discrete problem is solved using a Newton method.
After each Newton step, the correction is interpolated from the vertices to the quadrature points.
This is easily possible, since the corrections are elements of a tangent space (and hence a linear space).
A similar approach is used in~\cite{dornisch_klinkel_simeon:2013,dornisch_klinkel:2014} in the context of isogeometric analysis,
where NURBS basis functions are employed for a geometrically exact representation of the director vector at the quadrature points.
However, for a related model~\cite{simo_vu-quoc:1986}, \citet{crisfield_jelenic:1999} showed that
this approach leads to an artificial path dependence of the solution.  An additional disadvantage is
that discretization and solution algorithm are not clearly separated.  This makes analyzing the method
difficult.

One last approach regards the manifold $\text{SO}(3)$ as a submanifold
of a linear space.  One can then interpolate in this space, and project the result back onto the
manifold.  To the knowledge of the authors this approach has never been used for shell models.
For harmonic maps into the unit sphere it has been proposed and analyzed in \cite{bartels_prohl:2007}.
The approach is attractive for its simplicity.  However, the result of the discrete problem depends
on the embedding.  This is of particular importance in the case of rotations, which can be interpreted
as a submanifold of $\R^{3 \times 3}$ (in which case the projection is the polar decomposition),
but also (as quaternions) as a submanifold of $\R^4$ (see Section~\ref{sec:quaternion_coordinates}).
Furthermore, the approach has only been investigated for discretizations of first order,
and it is unclear whether higher approximation orders are possible as well.

\bigskip

In this article we propose a new discretization based on Geodesic Finite Elements, which solves most of the
shortcomings of the previous methods.  Geodesic Finite Elements (GFE), originally introduced in~\cite{sander:2012,sander:2013},
are a natural generalization of standard Lagrangian finite elements to spaces of functions mapping into a
general Riemannian manifold $M$.  The core idea is to write Lagrangian interpolation $T_\text{ref} \to \R$
of values $v_1,\dots,v_m \in \R$ on a reference element $T_\text{ref}$ as a minimization problem
\begin{equation*}
 \xi \mapsto \argmin_{w \in \R} \sum_{i=1}^m \lambda_i(\xi) \abs{v_i -w}^2,
\end{equation*}
where the $\lambda_1, \dots, \lambda_m : T_\text{ref} \to \R$ are the Lagrangian shape functions.
For values $v_1, \dots, v_m$ in a Riemannian manifold $M$, this formulation can be generalized using the
Riemannian distance
\begin{equation*}
 \xi \mapsto \argmin_{w \in M} \sum_{i=1}^m \lambda_i(\xi) \dist(v_i, w)^2.
\end{equation*}
This construction is also known as the Karcher mean~\cite{karcher:1977} or the Riemannian center of mass.
It forms the basis of a general finite element theory for functions mapping into a manifold $M$~\cite{sander:2012,sander:2013}.
Finite element spaces constructed this way are conforming in the sense that finite element functions belong
to the Sobolev space $W^{1,q}(\omega,M)$ for all $q \ge 2$.  Since their formulation is based on metric
properties of $M$, they are naturally equivariant under isometries of $M$.
Optimal a priori discretization error bounds have been given in~\cite{grohs_hardering_sander:2014}.

When using this technique for the case $M = \text{SO(3)}$ considered here (but the same holds also when
discretizing one-director models with $M=S^2$ such as the one proposed in~\cite{simo_fox:1989}), the resulting
discrete model has many desirable properties.  Since the FE spaces are conforming, there is no consistency
error introduced when evaluating the continuous energy for finite element functions.  Since no angles
and no ``special orientations'' appear in the discretization, the discrete model is not restricted to small
or moderate rotations.  Indeed, as we demonstrate in Section~\ref{sec:twisted_strip}, arbitrary
rotations in the deformation can be handled with ease.  Finally, from the equivariance of the nonlinear
interpolation follows that the frame invariance of the continuous model~\eqref{eq:general_energy}
is preserved by the discretization, and we obtain a completely frame-invariant discrete problem.

As an additional advantage, the fact that the FE space is contained in the continuous ansatz space
$H^1(\omega,\R^3) \times W^{1,q}(\omega,\text{SO(3)})$ implies that properties of the tangent matrix
can be inferred from corresponding properties of the continuous tangent operator.  In particular,
we directly obtain symmetry of the tangent matrix.  The tangent matrix is positive definite if the
continuous tangent operator is.

The algebraic formulation corresponding to the discrete problem is a minimization problem posed in the
product space $\mathcal{M} \colonequals \R^{3N} \times \SOdrei^N$, where $N$ is the number of
Lagrange nodes of the grid.  The space $\mathcal{M}$ is a $6N$-dimensional Riemannian manifold.
To solve this minimization problem we use a Riemannian trust-region algorithm~\cite{absil_mahony_sepulchre:2008},
which is a globalized Newton method.  As such, it is guaranteed to converge to at least a stationary point
of the algebraic energy for any initial iterate, and without using intermediate loading steps.
At each step of the method, a constrained quadratic minimization problem needs to be solved.
We propose to use a monotone multigrid method~\cite{kornhuber:1997,sander:2012}, which allows efficient and robust
solutions of the constrained problems even on fine grids.
As a variant of the Newton method, the trust-region algorithm requires tangent matrices
of the energy.  We obtain those matrices completely automatically by using automatic
differentiation (AD) as implemented in the software ADOL-C~\cite{walther_griewank:2012,griewank_walther:2008}.

\bigskip

In this article we show three numerical examples.  First we compute the post-critical behavior of an
$L$-shaped beam.  This was posed as a benchmark problem in~\cite{wriggers_gruttmann:1993,argyris_et_al:1979,simo_vu-quoc:1986,simo_fox_rifai:1990}, and we
compare our results with results given there.  Secondly, we demonstrate that our discretization
does indeed allow unrestricted rotations.  For this we simulate a long elastic strip, which we clamp
on one short end, and subject it to several full rotations at the other end.  Finally, to show that the
Cosserat shell model can represent non-classical microstructure effects, we use it to produce wrinkles
in a sheared rectangular membrane.  Such shearing tests have been performed experimentally by
\cite{wong_pellegrino:2006a}, and we obtain excellent quantitative agreement with their results.

\bigskip

This article is structured as follows:  In Chapter~\ref{sec:continuous_model} we present the continuous model
and discuss a few of its properties.  Chapter~\ref{sec:GFE} introduces the geodesic finite element method,
specialized for the case $M=\SOdrei$ needed for the Cosserat shell model.  Chapter~\ref{sec:discrete_problem}
discusses the resulting discrete and algebraic models.
Chapter~\ref{sec:numerical_minimization} explains the Riemannian trust-region method used to find
energy minimizers without loading steps.
Chapter~\ref{sec:numerical_tests} gives the three numerical examples.  Finally, an appendix
collects various important facts about $\text{SO(3)}$ needed to implement the GFE method.

\section{The continuous Cosserat shell model}
\label{sec:continuous_model}

In this chapter we present the planar Cosserat shell model and discuss its features.  The detailed derivation of the shell model from a three-dimensional
parent Cosserat model was presented in the papers \cite{Neff_plate04_cmt,Neff_plate07_m3as}.
The intermediate shell model for
infinitesimal strain is described in Section~\ref{sec:infinitesimal_shell_model}.  The complete finite-strain model
is then introduced in Section~\ref{sec:finite_strain_shell_model}.

\subsection{The small-strain planar Cosserat shell model}
\label{sec:infinitesimal_shell_model}

We consider a thin domain $\Omega_h \subset \R^3$ of the form $\Omega_h=\omega\times[-h/2 , h/2]$, where $\omega$ is a bounded domain in $\mathbb{R}^2$ with smooth boundary $\partial \omega$, and $h>0$ is the thickness of the planar shell.
The domain $\Omega_h$ is the region occupied by the reference configuration of the parent 3D Cosserat continuum.
Let $\{e_1,e_2,e_3\}$ be the unit vectors along the axes of the reference Cartesian coordinate system, denote by $\varphi:\Omega_h \to \R^3$ the deformation, and by $\overline{R}:\Omega_h\to\mathrm{SO}(3)$ the independent microrotation of this micropolar continuum.

For the planar shell model we want to find a reasonable approximation $(\varphi_s,\overline{R}_s)$ of $(\varphi,\overline{R})$  involving only two-dimensional quantities, i.e., expressed with the help of functions of the in-plane coordinates $(x,y)$. Therefore, we assume a quadratic ansatz in the thickness coordinate $z$ for the finite deformation $\varphi_s:\Omega_h\rightarrow  \mathbb{R}^3$
\begin{equation}\label{eq:quadratic_ansatz}
    \varphi_s(x,y,z)=m(x,y)+\Big(z\varrho_m(x,y)+\dfrac{z^2}{2}\,\varrho_b(x,y)\Big)\, \boldsymbol{d}(x,y).
\end{equation}
Here $m:\omega\to\mathbb{R}^3$ describes   the deformation of the midsurface of the shell,  and $\boldsymbol{d}:\omega\to\mathbb{R}^3$ is an independent unit director.
We assume the rotations $\overline{R}_s:\Omega_h\rightarrow\mathrm{SO}(3)$ for thin and homogeneous shells
to be independent of the thickness variable $z$, i.e.,
\begin{equation*}
    \overline{R}_s(x,y,z)=\overline{R}_s(x,y,0)\qquad\text{for}\quad z\in\big[-\frac{h}{2}\,,\,\frac{h}{2}\,\big],
\end{equation*}
and we specialize  the independent unit director $\boldsymbol{d}$ in the ansatz \eqref{eq:quadratic_ansatz} by choosing
\begin{equation*}
    \boldsymbol{d}(x,y)\colonequals\overline{R}_s(x,y,0) e_3 \equalscolon \overline{R}_3\,.
\end{equation*}
Thus, the director $\boldsymbol{d}(x,y)$ is taken as the third column of the orthogonal matrix $\overline{R}_s(x,y)$, and the model
now also includes drilling rotations about the director $\boldsymbol{d}$.  The drilling rotations  are determined by the first two columns of $\overline{R}_s$. For the sake of simplicity, we drop the index $s$  and write $\overline{R} $ instead of $\overline{R}_s$ in what follows.

When the director $\boldsymbol{d}(x,y)$ is not normal to the midsurface $m(x,y)$, then transverse shear deformation occurs. The scalar functions $\rho_m,\rho_b:\omega\rightarrow\mathbb{R}$ in \eqref{eq:quadratic_ansatz} describe the symmetric thickness stretch (for $\rho_m\neq 1$) and the asymmetric thickness stretch (for $\rho_b\neq 0$) about the midsurface. The scalar field $\rho_m$ is mainly membrane related, while $\rho_b$ is mainly bending related.  The fields have the following expressions \cite{Neff_plate04_cmt}
\begin{equation*}
    \begin{aligned}
      \varrho_m & \colonequals 1-\dfrac{\lambda}{2\mu+\lambda}\,\big[\langle\,(\nabla m|\,0),\overline{R}\,\rangle-2\big]+\dfrac{\langle\, N_{\mathrm{diff}}\,,\, \overline{R}_3 \,\rangle}{2\mu+\lambda}\,,
      \\
      \varrho_b & \colonequals -\dfrac{\lambda}{2\mu+\lambda}\,\langle\,(\nabla \overline{R}_3|\,0),\overline{R}\,\rangle+\dfrac{\langle\, N_{\mathrm{res}}\,,\, \overline{R}_3 \,\rangle}{(2\mu+\lambda)h}\,,
    \end{aligned}
\end{equation*}
where the parameters $\lambda,\mu > 0$ are the Lam\'e constants of classical isotropic elasticity, and $N_{\mathrm{res}},\,N_{\mathrm{diff}}:\omega\rightarrow\R^3$ are defined in terms of the prescribed tractions $N^{\mathrm{trans}}$ on the transverse boundaries $z=\pm h/2\,\,$ by
\begin{equation*}
    \begin{aligned}
       N_{\mathrm{res}}(x,y) & \colonequals\big[ N^{\mathrm{trans}}(x,y,\,\frac{h}{2}\,)+N^{\mathrm{trans}}(x,y,-\frac{h}{2}\,)\big],
       \\
       N_{\mathrm{diff}}(x,y) & \colonequals\frac{1}{2}\,\big[ N^{\mathrm{trans}}(x,y,\,\frac{h}{2}\,)-N^{\mathrm{trans}}(x,y,-\frac{h}{2}\,)\big].
    \end{aligned}
\end{equation*}

The strain measures for the planar Cosserat shell model are the following: the micropolar non-symmetric stretch tensor $\overline{U}$ is defined as
\begin{equation*}
\overline{U}
=
\overline{R}^T\hat{F}
\qquad\text{with}\qquad
\hat{F}=(\nabla m|\, \overline{R}_3)\in\mathbb{M}^{3\times3},
\end{equation*}
while the micropolar curvature tensor $\mathfrak{K}_s$  (of third order) and the micropolar bending tensor $\mathfrak{K}_b$ (of second order) are given by
\begin{align*}
      \mathfrak{K}_s & \colonequals \big( \,\overline{R}^T(\nabla \overline{R}_1|\,0)\,,\, \overline{R}^T(\nabla \overline{R}_2|\,0)\,,\, \overline{R}^T(\nabla \overline{R}_3|\,0)\,\big)
      \equalscolon
      \big( \,\mathfrak{K}_s^1\,,\, \mathfrak{K}_s^2\,,\, \mathfrak{K}_s^3\,\big)
      \in \mathbb{M}^{3 \times 3 \times 3},\\
      \mathfrak{K}_b & \colonequals  \overline{R}^T(\nabla \overline{R}_3|\,0)= \mathfrak{K}_s^3
      \in \mathbb{M}^{3 \times 3}\,.
\end{align*}
We have used the superposed caret and bars for $\hat{F}$, $\overline{R}$, $\overline{U}$ in order to distinguish these tensors from the classical notations in 3D elasticity for deformation gradient $F$, the continuum rotation $R=\operatorname{polar}(F)$, and the symmetric continuum stretch tensor $U=R^TF=\sqrt{F^TF}$.

We mention that the kinematical structure of this Cosserat shell model is in fact equivalent
to the kinematical structure of nonlinear 6-parameter resultant shell theory \cite{Pietraszkiewicz-book04,Libai98,Eremeyev06}, as it was pointed out in \cite{,Birsan-Neff-JElast-2013,Birsan-Neff-MMS-2014,Birsan-Neff-Danzig13,Birsan-Neff-L54-2014}.

As a result of the dimensional reduction procedure, the following two-dimensional minimization
problem for the deformation of the midsurface $m:\omega\to\R^3$ and the microrotation field
$\overline{R}:\omega \to \SOdrei$ is obtained~\cite{Neff_plate04_cmt}:
\begin{problem}
Find a pair $(m,\overline{R})$ which minimizes the functional
\begin{equation}\label{eq:small_strain_energy}
       I(m,\overline{R})=\int_\omega h\,W_{\mathrm{mp}}(\overline{U})+ h\,W_{\mathrm{curv}}(\mathfrak{K}_s)+ \dfrac{h^3}{12}\, W_{\mathrm{bend}}(\mathfrak{K}_b)\,\mathrm{d}\omega -  \Pi(m,\overline{R}_3)\,,
\end{equation}
subject to suitable boundary conditions for the deformation and rotation.
\end{problem}

The three parts of the total elastically stored energy density of the shell correspond to membrane-strain $W_{\text{mp}}$,
total  curvature-strain $W_{\text{curv}}$ and specific  bending-strain $W_{\text{bend}}$.
They have the expressions
\begin{align}
  \nonumber
  W_{\mathrm{mp}}(\overline{U}) & = \mu\|\sym(\overline{U}-1\!\!1)\|^2 +\mu_c\|\,\operatorname{skew}(\overline{U}-1\!\!1)\|^2+ \dfrac{\mu\lambda}{2\mu+\lambda}\,\tr\, \big[\sym(\overline{U}-1\!\!1)\big]^2
  \vspace{3pt}\\
  \nonumber
  \qquad & =  \mu\underbrace{\|\sym( ( \overline{R}_1|\,\overline{R}_2)^T\nabla m- 1\!\!1_2  )\|^2}_{\text{shear-stretch energy}} \,\,+\,\mu_c\underbrace{\|\operatorname{skew}(( \overline{R}_1|\,\overline{R}_2)^T\nabla m)\|^2}_{\text{first order drill energy}}
  \vspace{3pt}\\
  \label{e10}
  & \quad+\dfrac{(\mu+\mu_c)}{2}\,\underbrace{\kappa\,\big(\langle \overline{R}_3\,,\,m_x\rangle^2 + \langle \overline{R}_3\,,\,m_y\rangle^2\big)}_{\text{classical transverse shear energy}} + \dfrac{\mu\lambda}{2\mu+\lambda}\,\underbrace{\tr\big[\sym(( \overline{R}_1|\,\overline{R}_2)^T\nabla m- 1\!\!1_2)\big]^2}_{\text{volumetric stretch energy}} \,,
  \vspace{6pt}\\
  \displaybreak[0]
  \nonumber
  W_{\mathrm{curv}}(\mathfrak{K}_s) & = \mu\, L_c^q\, \|\mathfrak{K}_s\,\|^q= \mu\, L_c^q\,\Big(\|\mathfrak{K}_s^1\,\|^2+\|\mathfrak{K}_s^2\,\|^2+ \|\mathfrak{K}_s^3\,\|^2\Big)^{q/2}\,,\vspace{3pt}\\
  \label{eq:bending_energy}
  W_{\mathrm{bend}}(\mathfrak{K}_b) & = \mu\, \|\sym(\mathfrak{K}_b)\|^2 +\mu_c\|\operatorname{skew}(\mathfrak{K}_b)\|^2+\,\dfrac{\mu\lambda}{2\mu+\lambda} \,\tr\big[\sym(\mathfrak{K}_b)\big]^2 ,
\end{align}
where the additional parameter
$\mu_c\geq 0 $ is called the Cosserat couple modulus, and $\kappa$ is a shear correction factor ($0<\kappa\leq 1$).  For $\mu_c > 0$ the elastic strain energy density $W_{\mathrm{mp}}(\overline{U})$ is
uniformly convex in $\overline{U}$, but for the important case $\mu_c = 0$ this property is lost. Therefore, the case $\mu_c = 0$ must be investigated separately.  In the curvature energy density $W_\text{curv}\,$, the parameter $L_c > 0$   is an internal length which is characteristic for the material, and is responsible for size effects. Note that this is a first-order model,
i.e., no second or higher derivatives of the independent variables $m$ and $\overline{R}$ appear.
Also, the energy depends on the midsurface deformation $m$ and microrotations $\overline{R}$
only through the frame-indifferent measures $\overline{U}$ and $\mathfrak{K}_s$\,. Thus, in the absence of external forces, the planar shell model is fully frame-indifferent in the sense that
\begin{equation*}
I(Qm,Q\overline{R}) = I(m,\overline{R}),\qquad \forall\,Q\in\SOdrei.
\end{equation*}

The reduced external loading functional $\Pi(m,\overline{R}_3)$ appearing in \eqref{eq:small_strain_energy} is a linear form in
$(m,\overline{R}_3)$, defined in terms of the underlying three-dimensional loads by
\begin{equation*}
\Pi(m,\overline{R}_3)=\int_\omega\langle\,\overline{f},m\rangle+ \langle\,\overline{M},\overline{R}_3\rangle \,\mathrm{d}\omega+\int_{\gamma_s} \langle\,\overline{N},m\rangle+\langle\,\overline{M}_c,\overline{R}_3\rangle \,ds,
\end{equation*}
where $\gamma_s\times [-\frac{h}{2}\,,\,\frac{h}{2}]\subset\partial\omega\times[-\frac{h}{2}\,,\,\frac{h}{2}]$ is the part of the lateral boundary of $\Omega_h$ where external surface forces and couples are prescribed. The vector fields $\overline{f}, \overline{M}, \overline{N} $ and $\overline{M}_c$ denote the resultant body force, resultant body couple, resultant surface traction and resultant surface couple, respectively~\cite{Neff_plate04_cmt}.

For the Dirichlet boundary conditions we suppose that there exists a prescribed function
$g_d : \Omega_h \to \R^3$, whose restriction to the Dirichlet part of the boundary gives
the prescribed displacement.  We further introduce the abbreviation
\begin{equation*}
 g'_d : \omega \to \R^3,
 \qquad
 g'_d(x,y) \colonequals \nabla g_d(x,y,0) e_3.
\end{equation*}

For the midsurface deformation $m$ we then consider the boundary conditions
\begin{equation}\label{e11}
    m(x,y)_{\big|\gamma_0}=g_d(x,y,0),
\end{equation}
on the Dirichlet part $\gamma_0$ of the boundary $\partial \omega$.

For the  microrotations $\overline{R}$ we can consider various possible alternative boundary conditions on $\gamma_0\,$,
see \cite{Neff_plate04_cmt,Neff_plate07_m3as}. In what follows, we consider two types:
\begin{align}
  \label{e12}
  & \text{1. \emph{free boundary conditions on} } \overline{R}  ,\quad \text{i.e., induced Neumann type (natural) conditions};
  \vspace{3pt}\\
  \label{e14}
  & \text{2. \emph{rigid director prescription},\quad i.e.,} \quad {\overline{R}_3}_{\big|\gamma_0}=\frac{g'_d}{\|\,g'_d\,\|},\\
  \nonumber
  & \qquad \text{together with zero Neumann conditions for the drilling degree of freedom}.
 \end{align}

The existence of minimizers for this Cosserat planar shell model under various assumptions on the coefficients and boundary conditions has been proved in \cite{Neff_plate04_cmt,Neff_plate07_m3as}. For instance, in the case when the Cosserat couple modulus  is positive ($\mu_c>0$) and for rigid director prescription boundary conditions \eqref{e14} on $\gamma_0$, the following existence result has
been shown in \cite{Neff_plate04_cmt}, using the direct method of the calculus of variations.

\begin{theorem}\label{thm:existence_infinitesimal_strain}
Let $\omega\subset\mathbb{R}^2$ be a bounded Lipschitz
domain, and assume that the material parameters satisfy
\begin{equation*}
\mu_c>0,\qquad   q\geq 2\,.
\end{equation*}
Let the  boundary data and external loads functions satisfy the regularity conditions
\begin{equation}\label{e15}
      g_d(x,y,0)\in H^1(\omega,\mathbb{R}^3),\qquad \operatorname{polar}\big(\nabla g_d(x,y,0)\big)\in W^{1,q}(\omega, \SOdrei),
\end{equation}
\begin{equation*}
      \overline{f}\in L^2(\omega,\mathbb{R}^3),\quad \overline{M}\in L^1(\omega,\mathbb{R}^3),\quad \overline{N} \in L^2(\gamma_s,\mathbb{R}^3) ,\quad \overline{M}_c\in L^1(\gamma_s,\mathbb{R}^3).
\end{equation*}
Then the minimization  problem \eqref{eq:small_strain_energy}--\eqref{eq:bending_energy}  with boundary conditions \eqref{e11} and \eqref{e14} admits at least one minimizing solution pair $(m,\overline{R})\in H^1(\omega,\mathbb{R}^3) \times W^{1,q}(\omega, \mathrm{SO}(3))$.
\end{theorem}

In the case of zero Cosserat couple modulus   ($\mu_c=0$) the mathematical treatment of the minimization problem is more difficult, due to the lack of unqualified coercivity  of the energy function with respect to the midsurface deformation $m$. The corresponding existence result for this case has been proved in \cite{Neff_plate07_m3as} using a new extended Korn's first inequality for plates and elasto-plastic shells \cite{Neff-Korn-02,Pompe-03}.
In this case, we need $q$ to be strictly larger than $2$.  However, the numerical evidence in Chapter~\ref{sec:numerical_tests}
suggests that existence also holds for $q=2$.
For the sake of simplicity, we present this result in the case of zero external loads, i.e., $\overline{f}=0$, $\overline{M}=0$, $\overline{N} =0$, $\overline{M}_c=0$.

\begin{theorem}\label{thm:existence_infinitesimal_strain_II}
Let $\omega\subset\mathbb{R}^2$ be a bounded Lipschitz
domain and assume that the material parameters satisfy
\begin{equation*}
\mu_c=0,\qquad   q> 2\,.
\end{equation*}
Let the  boundary data satisfy the regularity conditions
\begin{equation*}
      g_d(x,y,0)\in H^1(\omega,\mathbb{R}^3),\qquad \operatorname{polar}\big(\nabla g_d(x,y,0)\big)\in W^{1,q}(\omega, \SOdrei).
\end{equation*}
Then the minimization problem for the functional \eqref{eq:small_strain_energy}--\eqref{eq:bending_energy}  with boundary conditions \eqref{e11} and \eqref{e14} admits at least one minimizing solution pair $(m,\overline{R})\in H^1(\omega,\mathbb{R}^3) \times W^{1,q}(\omega, \mathrm{SO}(3))$.
\end{theorem}
The statement of Theorem~\ref{thm:existence_infinitesimal_strain_II} holds also in the case of non-vanishing external loads. In this respect, see the paper \cite{Neff_plate07_m3as}, where a modification of the external loading potential has been used.

\bigskip

Of particular interest is the choice of the new material parameters $\mu_c$ (the Cosserat couple modulus) and $L_c$.
Our model is derived from a 3D-Cosserat model in which the Cosserat couple modulus appears traditionally.
It controls the skew-symmetric part of the stresses, and enforces $\overline{R} = \operatorname{polar}(\hat{F})$
for the limit case $\mu_c \to \infty$.
From the literature, there does not exist a single material for which the value of the parameter $\mu_c$ has been identified unambiguously.
Considering this situation, in~\cite{Neff_zamm06} it is argued that this parameter must be set to zero when
modeling a continuous body. In~\cite{Neff_Ghiba_dynamic_micromorph13,Neff_Ghiba_unified_micromorph13}
the same question has been discussed in the larger framework of (infinitesimal) micromorphic continua with the same result:
the absence of $\mu_c$ leads to a more stringent physical description. Indeed, it implies that a linear Cosserat model collapses into classical linear elasticity.

However, in a geometrically nonlinear context, which is our case, a vanishing Cosserat couple modulus only implies that there is no
first-order coupling between rotations and deformation gradients \cite{Neff_Biot07}.
Compared with the classical Reissner--Mindlin kinematics without drill energy~\cite{Neff_Hong_Reissner08},
setting $\mu_c=0$ appears again as the most plausible choice. Since, therefore, there is no specific reason to have $\mu_c>0$, we omit this parameter.

The internal length $L_c$ appears in Cosserat models as a measure of the length scale of the material microstructure.
The numerical results of Chapter~\ref{sec:numerical_tests} show that values of $L_c$ in the micrometer range lead to realistic results.
However, we also note that the shell model with
$L_c\gg h$ can be useful for the description of graphene-sheets which have practically zero thickness but still
show a bending stiffness. In a classical shell model, we would expect zero bending resistance.

\subsection{A modified large strain Cosserat shell model}
\label{sec:finite_strain_shell_model}

We observe that the planar shell model presented above is appropriate for finite rotations, but only small elastic membrane strains,
since the membrane part $W_{\mathrm{mp}}$ of the energy density $I$ is quadratic. We now slightly  generalize the model to allow for
large elastic stretch as well. We consider again a minimization problem for the
energy functional
\begin{equation}
 \label{eq:finite_strain_energy}
 I(m,\overline{R})
 =
 \int_\omega h\,W_{\mathrm{mp}}(\overline{U})
           + h\,W_{\mathrm{curv}}(\mathfrak{K}_s)
           + \dfrac{h^3}{12}\, W_{\mathrm{bend}}(\mathfrak{K}_b)\,\mathrm{d}\omega
           -  \Pi(m,\overline{R}_3)\,,
\end{equation}
but we replace the membrane part of $I$ by
\begin{align}
 \nonumber
      W_{\mathrm{mp}}(\overline{U}) & = \mu\|\sym(\overline{U}-1\!\!1)\|^2 +\mu_c\|\operatorname{skew}(\overline{U}-1\!\!1)\|^2+ \dfrac{\mu\lambda}{2\mu+\lambda}\,\dfrac{1}{2}\,\Big((\det \overline{U}-1 )^2+ \big((\det \overline{U})^{-1}- 1 \big)^2\Big) \vspace{3pt}\\
       \nonumber
      \qquad & =  \mu\underbrace{\|\sym( ( \overline{R}_1|\,\overline{R}_2)^T\nabla m- 1\!\!1_2  )\|^2}_{\text{shear-stretch energy}} \,\,+\,\mu_c\underbrace{\|\operatorname{skew}(( \overline{R}_1|\,\overline{R}_2)^T\nabla m)\|^2}_{\text{first order drill energy}}
  \vspace{3pt}\\ 
  \label{eq:finite_strain_membrane_energy}
  & \quad+\dfrac{(\mu+\mu_c)}{2}\,\underbrace{\kappa\,\big(\langle \overline{R}_3\,,\,m_x\rangle^2 + \langle \overline{R}_3\,,\,m_y\rangle^2\big)}_{\text{classical transverse shear energy}}
   \vspace{3pt}\\
  \nonumber
  & \quad+\dfrac{\mu\lambda}{2\mu+\lambda}\,\underbrace{\dfrac{1}{2}\,\Big(\big(\det(\nabla m\,|\overline{R}_3) -1 \big)^2+ \big(\det(\nabla m\,|\overline{R}_3)^{-1}- 1 \big)^2\Big)}_{\text{modified volumetric stretch response}} \,.
\end{align}
In this expression, we have replaced the quadratic volumetric stretch part $\tr[\sym(\overline{U}-\mathbbm{1})]^2$
of~\eqref{e10} by the non-quadratic expression
\begin{equation*}
 \frac{1}{2}\,\Big((\det \overline{U}-1 )^2+ \big((\det \overline{U})^{-1}- 1 \big)^2\Big),
\end{equation*}
which is volumetrically exact.
However, since
\begin{equation*}
\frac{1}{2}\,\Big((\det \overline{U}-1 )^2+ \big( (\det \overline{U})^{-1}- 1 \big)^2\Big)
=
\tr \big[ \sym(\overline{U}-\mathbbm{1})\big]^2 +O\big(\| \overline{U}-\mathbbm{1}\|^3\big),
\end{equation*}
the quadratic membrane energy \eqref{e10} of the previous section can be recovered by linearization
at $\mathbbm{1} \in \mathbb{M}^{3 \times 3}$.

For the nonlinear modified model \eqref{eq:finite_strain_energy} we set the following expression  for the modified thickness stretch
\begin{equation*}
\varrho_m
\colonequals
\dfrac{1}{1+\frac{\lambda}{2\mu+\lambda}\,(\det \overline{U}-1)}\,\in (0,\infty),
\end{equation*}
which can be used for the a posteriori reconstruction of the bulk deformation.

The modified membrane energy density \eqref{eq:finite_strain_membrane_energy}  represents an improvement over the initial planar shell model \eqref{e10} in various regards. Indeed, we note that
\begin{equation*}
W_{\mathrm{mp}}(\overline{U})\rightarrow\infty
\qquad\text{if}\quad
\det \overline{U}\rightarrow 0.
\end{equation*}
Moreover, for any fixed $\overline{R}$ the energy $W_{\mathrm{mp}}$  is polyconvex \cite{Schroeder_Neff_Ebbing07,Ebbing_Schroeder_Neff_AAM09,Ebbing_Balzani_Schroeder_Neff_CMS09,Balzani_Schroeder_Neff_IJSS06,Neff_Schroeder_CISM_07} with respect to $\nabla m$,
and it is uniformly Legendre--Hadamard elliptic, independently of $\mu_c\geq 0$.

The following existence result for the modified model, in the important case $\mu_c = 0$, was originally proved in \cite{Neff_plate07_m3as}. Again, we assume vanishing external loads for simplicity.

\begin{theorem}\label{thm:existence_finite_strain}
Let $\omega\subset\mathbb{R}^2$ be a  bounded Lipschitz
domain and assume that the boundary data  satisfies~\eqref{e15}. 

Then the minimization  problem for the functional \eqref{eq:finite_strain_energy} with the parameters
$$\mu_c=0\qquad and \qquad q>2,$$
with boundary conditions \eqref{e11}, \eqref{e14}  admits at least one minimizing solution pair $(m,\overline{R})\in H^1(\omega,\mathbb{R}^3) \times W^{1,q}(\omega, \SOdrei)$, with
$$\,\det \big(\nabla m\,|\,\overline{R}_3\big)=\det\,\hat{F}>0\,$$ 
almost everywhere in $\omega$.
\end{theorem}

We note that the formulation \eqref{eq:finite_strain_energy} has the same linearized behavior as the initial model \eqref{eq:small_strain_energy} and it reduces upon linearization to the classical infinitesimal-displacement Reissner--Mindlin model for the choice of parameters $\,\mu_c=0\,$ and $\,q>2$.

\begin{remark}
The Cosserat model presented above can be extended to a general nonplanar shell model.
Indeed, instead of the domain $\Omega_h$ and the ansatz for plates \eqref{eq:quadratic_ansatz}, one can begin with a shell-like (curved) thin domain and an appropriate ansatz for shells. Then, the formal dimensional reduction to a two-dimensional shell model is derived analogously as in the case of plates, but involves additional tools from classical differential geometry of surfaces for the description of shell configurations. The resulting Cosserat shell model is quite general and has the advantage that it can be used to also describe elasto-plastic and visco-elasto-plastic material behavior. This work is currently in progress.
\end{remark}

\section{Geodesic finite elements}\label{sec:GFE}

Discretization of the shell models presented in the previous section is difficult, because the orientation configuration space
$W^{1,q}(\omega,\text{SO(3)})$ is not linear.  As a consequence, linear,
and more generally polynomial, interpolation is undefined in these spaces, and standard finite element
methods cannot be used.

Geodesic finite elements (GFE) are a generalization of standard finite elements to problems for functions
with values in a nonlinear Riemannian manifold $M$.  We give a brief introduction and state the relevant
features without proof.  While geodesic finite elements can be constructed easily for very general $M$,
we state all results here for the case $M = \text{SO(3)}$ only.
The interested reader is referred to the original publications \cite{sander:2012,sander:2013}
for more details.

The definition of GFE spaces consists of two parts.  First, nonlinear interpolation
functions are constructed that interpolate values given on a reference element.
Then, these interpolation functions are pieced together to form global finite element
spaces for a given grid.

\subsection{Geodesic interpolation}

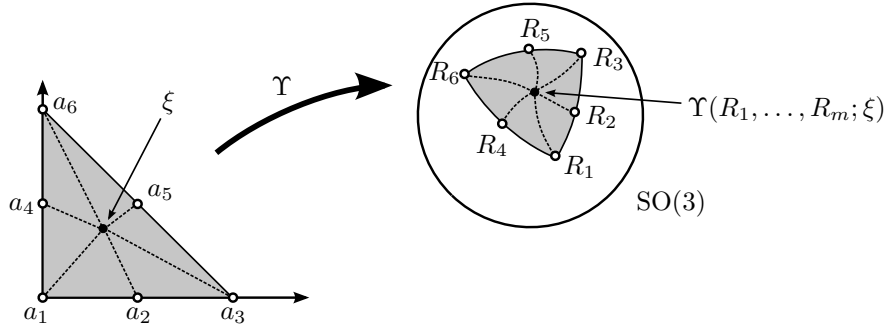
\begin{figure}
\begin{center}
 \begin{overpic}[width=0.6\textwidth]{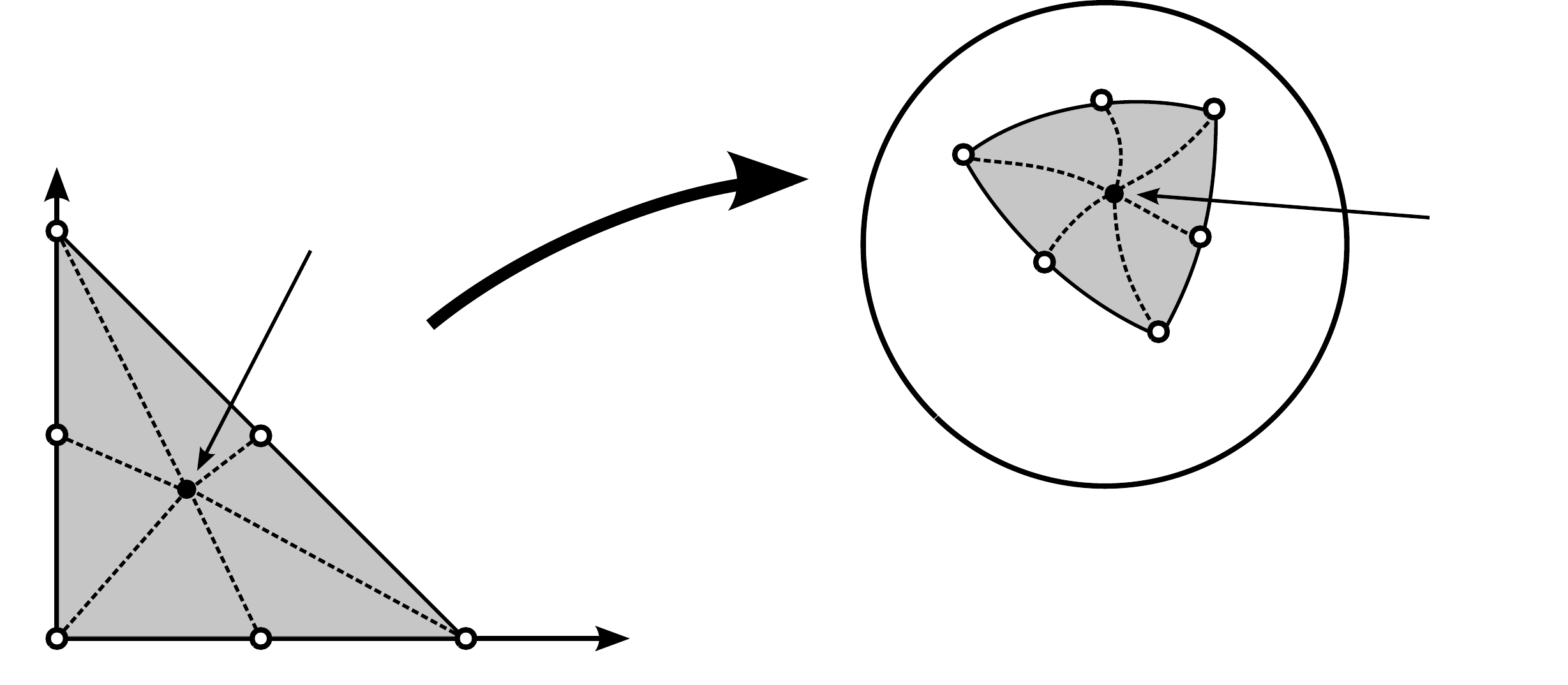}
  \put( 1, 0){$a_1$}
  \put(15, 0){$a_2$}
  \put(28, 0){$a_3$}
  \put(-1,15){$a_4$}
  \put(18,17){$a_5$}
  \put( 5,29){$a_6$}
  \put(75,20){$R_1$}
  \put(78,27){$R_2$}
  \put(79,35){$R_3$}
  \put(63,23){$R_4$}
  \put(69,39){$R_5$}
  \put(56.5,33){$R_6$}
  \put(20,29){$\xi$}
  \put(35,31){$\Upsilon$}
  \put(85,15){$\SOdrei$}
  \put(92,28){$\Upsilon(R_1,\dots,R_m;\xi)$}
 \end{overpic}
\end{center}
\caption{Second-order geodesic interpolation from the reference triangle into a sphere}
 \label{fig:geodesic_interpolation}
\end{figure}

We focus on the case of a two-dimensional domain $\omega$.  All constructions and results work
mutatis mutandis also for domains of other dimensions.

Let $T_\text{ref}$ be a triangle or quadrilateral in $\R^2$.  We call $T_\text{ref}$ the reference
element.  On $T_\text{ref}$ we assume the existence of a set of $p$-th order Lagrangian interpolation polynomials,
i.e., a set of Lagrange nodes $a_i \in T_\text{ref}$, $i=1,\dots,m$, and corresponding
polynomial functions $\lambda_i : T_\text{ref} \to \R$ of order $p$ such that
\begin{equation*}
 \lambda_i(a_j) = \delta_{ij}\quad \text{for $i,j=1,\dots,m$},
 \qquad \text{and} \qquad
 \sum_{i=1}^m \lambda_i \equiv 1.
\end{equation*}
We want to generalize Lagrangian interpolation to the case of values $R_1,\dots,R_m \in \SOdrei$ associated
to the Lagrange nodes $a_i$.  In other words, we want to construct a function
$\Upsilon : T_\text{ref} \to \SOdrei$ such that $\Upsilon(a_i) = R_i$ for all $i = 1,\dots,m$.
This is a non-trivial task because $\SOdrei$ is not a vector space.

To motivate our construction we note that the usual Lagrangian interpolation of values $v_1,\dots,v_m$ in $\R$ can be written as a minimization
problem
\begin{equation*}
 \xi \mapsto \argmin_{w \in \R} \sum_{i=1}^m \lambda_i(\xi) \abs{v_i -w}^2
\end{equation*}
for each $\xi \in T_\text{ref}$.
This formulation can be generalized to values in $\SOdrei$.  We use $\dist(\cdot,\cdot)$
for the canonical (geodesic) distance on $\SOdrei$, which is
\begin{equation*}
 \dist(R_1,R_2)
 =
 \norm{\log R_1^T R_2}.
\end{equation*}

\begin{definition}[\cite{sander:2013}]
\label{def:geodesic_interpolation}
 Let $\{\lambda_i\}_{i=1}^m$ be a set of  $p$-th order
 scalar Lagrangian shape functions on the reference element $T_\text{ref}$, and let $R_i \in \SOdrei$,
 $i=1,\dots,m$ be values at the corresponding Lagrange nodes.  We call
\begin{align*}
 \Upsilon & \; : \; \SOdrei^m \times T_\text{ref} \to \SOdrei \\
 \Upsilon(R_1,\dots,R_m;\xi) & = \argmin_{R \in \SOdrei} \sum_{i=1}^m \lambda_i(\xi) \dist(R_i,R)^2
\end{align*}
$p$-th order geodesic interpolation on $\SOdrei$.
\end{definition}
To make the construction easier to understand we work out a simple example.
\begin{example}
 Let $T_\text{ref}$ be the reference triangle
 \begin{equation*}
  T_\text{ref}
  \colonequals
  \big\{ \xi \in \R^2 \; : \; \xi_1 \ge 0, \; \xi_2 \ge 0, \; \xi_1 + \xi_2 \le 1\big\},
 \end{equation*}
 and consider the first-order case $p=1$.  In this case, the Lagrange nodes $a_1$, $a_2$, $a_3$ are the triangle
 vertices, and the corresponding shape functions are
 \begin{equation*}
  \lambda_1(\xi) = 1 - \xi_1 -\xi_2,
  \qquad
  \lambda_2(\xi) = \xi_1,
  \qquad
  \lambda_2(\xi) = \xi_2.
 \end{equation*}
 These are simply the barycentric coordinates of $\xi$ with respect to $T_\text{ref}$.  Let $R_1, R_2, R_3$ be
 given values on $\SOdrei$.
 The image of $T_\text{ref}$ under $\Upsilon$ is then a (possibly degenerate) geodesic triangle on $\SOdrei$
 with corners $R_1, R_2, R_3$.  In particular, the edges of $T_\text{ref}$ map onto geodesics on $\SOdrei$ (\cite[Lem.\,2.2 with Cor.\,2.2]{sander:2012}).
 Even more, the map $\Upsilon$ is equivariant under permutations of the values $R_1, R_2, R_3$
 (\cite[Lem.\,4.3]{sander:2013}), a property not shared by various other commonly used discretization
 techniques~\cite{muench:2007,mueller:2009,simo_fox_rifai:1990}.  A visualization of this interpolation
 function can be found in~\cite{sander:2014}.  Also, Figure~\ref{fig:geodesic_interpolation} shows the corresponding
 second-order case.
\end{example}

\bigskip

While Definition~\ref{def:geodesic_interpolation} is an obvious generalization of Lagrangian interpolation
in linear spaces, it is by no means clear that it leads to a well-defined interpolation function 
for all coefficient sets $R_1,\dots, R_m \in \SOdrei$ and $\xi \in T_\text{ref}$.  Intuitively,
for fixed $\xi \in T_\text{ref}$, one would expect the functional
\begin{equation}
\label{eq:center_of_mass_energy}
  f_{\xi} : R \mapsto \sum_{i=1}^m \lambda_i(\xi) \dist(R_i,R)^2
\end{equation}
to have a unique minimizer if the $R_i \in \SOdrei$ are close enough to
each other in a certain sense. For the first-order case $p=1$, where all $\lambda_i$
are  non-negative on $T_\text{ref}$, this follows from a classic result of Karcher~\cite{karcher:1977},
which was later strengthened by Kendall~\cite{kendall:1990} (see also \cite{groisser:2004}).
Note that $\text{SO(3)}$ is complete and has constant sectional curvature of 1~\cite[Thm.\,2.7.1]{wolf:1974}.
\begin{theorem}[Kendall~\cite{kendall:1990}]
\label{thm:unique_minimizer_first_order}
Let $B_\rho$ be an open geodesic ball of radius $\rho< \pi / 2$ in $\text{SO(3)}$, and $R_1,\dots,R_m \in B_\rho$.
Let $\{ \lambda_i\}_{i=1}^m$ be a set of first-order Lagrangian shape functions.
Then the function
\begin{equation*}
  f_\xi : R \mapsto \sum_{i=1}^m \lambda_i(\xi) \dist(R_i,R)^2
\end{equation*}
has a unique minimizer in $B_\rho$ for all $\xi \in T_\text{ref}$.
\end{theorem}

If the polynomial order $p$ is larger than $1$, the weights $\lambda_i$ attain negative values on $T_\text{ref}$, and
the results of Karcher and Kendall cannot be used anymore.  Having all $R_i$ in a convex ball still guarantees
existence of a unique minimizer, but that minimizer may only be contained in a ball of larger size.
\begin{theorem}[\citet{sander:2013}]
\label{thm:well_posedness_vague}
 Let $B_D \subset B_\rho$
 be two concentric geodesic balls in $\SOdrei$ of radii $D$ and $\rho$, respectively,
 and let $R_1, \dots,R_m \in \SOdrei$.
 There are numbers $D$ and $\rho$ such that if $R_1,\dots,R_m \in B_D$, then the functional~\eqref{eq:center_of_mass_energy}
 has a unique minimizer in $B_\rho$.
\end{theorem}
A quantitative version of this result is given as Theorem~3.19 in~\cite{sander:2013}.  Unfortunately is is quite technical
and we have chose to omit it here.  When preparing the numerical examples of Chapter~\ref{sec:numerical_tests},
we have not encountered any problems stemming from a possible ill-posedness of the interpolation
for extreme configurations of the $R_1,\dots,R_m$.

\bigskip

To be able to use the interpolation functions as the basis of a finite element theory, they need to have
sufficient regularity.  The following result follows directly from the implicit function theorem.
\begin{theorem}
\label{thm:geodesic_interpolation_is_differentiable}
Let $R_1,\dots,R_m$ be coefficients on $\SOdrei$ with respect to
a $p$-th order Lagrange basis $\{ \lambda_i \}$ on a domain $T_\text{ref}$.
Under the assumptions of Theorem~\ref{thm:well_posedness_vague}, the function
\begin{equation*}
\Upsilon(R_1,\dots,R_m; \xi) \; : \; \SOdrei^m \times T_\text{ref} \to \SOdrei
\end{equation*}
 is infinitely differentiable with respect to the $R_i$ and $\xi$.
\end{theorem}
This result is proved in~\cite{sander:2012,sander:2013} for interpolation in general manifolds.

\subsection{Geodesic finite element functions}

The interpolation functions of the previous section can be used to construct
a generalization of Lagrangian finite element spaces to functions with values in $\SOdrei$.

For this, let $\omega$ be the two-dimensional parameter domain of our planar Cosserat shell model, and suppose it has
piecewise linear boundary.  Let $\mathcal{G}$
be a conforming grid for $\omega$ with triangle and/or quadrilateral elements.
Let $n_i \in \omega$, $i=1,\dots,N$ be a set of Lagrange nodes
such that for each element $T$ of $\mathcal{G}$ there are $m$ nodes
$a_{T,i}$ contained in $T$, and such that the $p$-th order interpolation problem
on $T$ is well posed.

\begin{definition}[Geodesic Finite Elements \cite{sander:2013}]
Let $\mathcal{G}$ be a conforming grid on $\omega$.  We call $R_h : \omega \to \SOdrei$
a geodesic finite element function if it is continuous, and
for each element $T \in \mathcal{G}$ the restriction $R_h|_T$
is a geodesic interpolation in the sense that
\begin{equation*}
 R_h|_T(x) = \Upsilon \big(R_{T,1}, \dots, R_{T,m}; \mathcal{F}_T(x) \big),
\end{equation*}
where $\mathcal{F}_T : T \to T_\text{ref}$ is affine or multilinear and the $R_{T,i}$ are values in $\SOdrei$
corresponding to the Lagrange nodes $a_{T,i}$.
  The space of all such functions $R_h$ will be denoted by $V_{p,h}^{\SOdrei}$.
\end{definition}

This construction has various desirable properties.
As a first result we note that the functions constructed in this way are $W^{1,q}$-conforming
for all $q \ge 2$.  This follows from a slight generalization of the proof for Theorem~3.1 in~\cite{sander:2012}.
\begin{theorem}
\label{thm:conformity}
 $V_{p,h}^{\SOdrei}(\omega) \subset W^{1,q}(\omega,\SOdrei)$ for all $p \ge 1$, $q \ge 2$.
\end{theorem}
Hence discrete approximation functions for the Cosserat microrotation field $\overline{R} : \omega \to \SOdrei$ are elements of the
space $W^{1,q}(\omega,\SOdrei)$, in which the Cosserat shell problem is well posed (Theorems~\ref{thm:existence_infinitesimal_strain}
and~\ref{thm:existence_finite_strain}).  This means that the energies~\eqref{eq:small_strain_energy}
and~\eqref{eq:finite_strain_energy} can be directly evaluated for geodesic finite element functions,
which simplifies the analysis considerably.

\bigskip

Since geodesic finite elements are defined using metric properties of $\SOdrei$ alone, we naturally
get the following equivariance result.
\begin{lemma}
\label{lem:equivariance_of_gfe}
 Let $O(3)$ be the orthogonal group on $\R^3$, which acts isometrically on $\SOdrei$ by left multiplication.
 Pick any element $Q \in O(3)$.  For any geodesic finite element function $R_h \in V_{p,h}^{\SOdrei}$ we define
 $QR_h : \omega \to \SOdrei$ by $(QR_h)(x) = Q(R_h(x))$ for all $x \in \omega$.  Then $QR_h \in V_{p,h}^{\SOdrei}$.
\end{lemma}
This lemma forms the basis of the frame-invariance of our discrete Cosserat shell model.

\bigskip

Optimal discretization error bounds for general GFE problems have been proved in~\cite{grohs_hardering_sander:2014}.
The application of those abstract results to the energy functionals considered in this paper
will be left for future work.

\section{Discrete and algebraic Cosserat planar shell problem}
\label{sec:discrete_problem}

We now discuss the minimization problem obtained by discretizing the continuous Cosserat shell model of Section~\ref{sec:continuous_model}
by geodesic finite elements.  For that, assume that the two-dimensional domain $\omega$ is discretized by a grid containing
triangle and/or quadrilateral elements.  For simplicity, we again assume that the domain boundary is resolved by the grid.
We also assume that the grid resolves the Dirichlet boundary~$\gamma_0$.

\subsection{The discrete problem}

The functional $I$ given in~\eqref{eq:finite_strain_energy} is
defined on the Cartesian product of the spaces $H^1(\omega,\R^3)$ and $W^{1,q}(\omega,\text{SO(3)})$.
The first factor is a standard Sobolev space of vector-valued functions.  For its discretization we introduce
the space $V_{p_1,h}^{\R^3}$ of conforming Lagrangian finite elements of $p_1$-th order with values in $\R^3$.
In the following we write $m_h$ for discrete displacement functions from $V_{p_1,h}^{\R^3}$.
For the rotation degree
of freedom $\overline{R} : \omega \to \text{SO(3)}$ we use the geodesic finite elements described in the previous chapter.  Denote by $V_{p_2,h}^\text{SO(3)}$
the $p_2$-th order GFE space for functions on $\omega$ with respect to the grid, and with values in $\text{SO(3)}$.
In the following we write $\overline{R}_h$ for discrete microrotations from $V_{p_2,h}^\text{SO(3)}$.

It is well known that $V_{p_1,h}^{\R^3} \subset H^1(\omega,\R^3)$ (see, e.g., \cite[Satz~5.2]{braess:2013}).  Additionally, we know
from Theorem~\ref{thm:conformity} that the finite element space  $V_{p_2,h}^\text{SO(3)}$ is a subset of $W^{1,q}(\omega,\text{SO(3)})$
for all $p_2 \in \mathbb{N}$.  Therefore, the energy functional $I$ is well defined on the product space
$\mathbf{V}_h \colonequals V_{p_1,h}^{\R^3} \times V_{p_2,h}^\text{SO(3)}$ for all $p_1,p_2 \in \mathbb{N}$.
A suitable discrete approximation of the geometrically nonlinear planar Cosserat shell model therefore consists
of the unmodified energy functional $I$ restricted to the space $\mathbf{V}_h$.

In analogy to the continuous model,
we consider the following boundary conditions for the discrete problem.  Let $g_{d,h} \in V_{h,p_1}^{\R^3}$ be a finite element approximation
of the Dirichlet boundary value function $g_d : \omega \to \R^3$, and let $g'_{d,h} \in V_{h,p_2}^{\R^3}$ be an
approximation of the vector field $g'_d$. Then we demand that the discrete displacement $m_h$ fulfill the condition
\begin{equation}
\label{eq:discrete_dirichlet_condition}
 m_h(x,y) = g_{d,h}(x,y)
 \qquad
 \text{for all $(x,y) \in \gamma_0$}.
\end{equation}
For the microrotations $\overline{R}$ we can define discrete approximations of the boundary conditions
\eqref{e12} and \eqref{e14}:  We either leave them free, corresponding to homogeneous Neumann conditions for $\overline{R}$,
or, alternatively, corresponding to \eqref{e14}, we can specify the direction of the transversal director vector $\overline{R}_3$
(rigid director prescription)
\begin{equation}
     \label{eq:discrete_rigid_director}
 (\overline{R}_h)_3 (x,y) = \frac{g'_{d,h}(x,y)}{\big\lVert g'_{d,h}(x,y)\big\rVert}
 \qquad
 \text{for all $(x,y) \in \gamma_0$}.
\end{equation}

Summing up, the discrete Cosserat shell problem is:
\begin{problem}[Discrete Cosserat shell problem]
Find a pair of functions $(m_h, \overline{R}_h)$ with
$m_h \in V_{p_1,h}^{\R^3}$ and $\overline{R}_h \in V_{p_2,h}^{\SOdrei}$ that minimizes the functional $I$ given in~\eqref{eq:finite_strain_energy},
subject to the constraints~\eqref{eq:discrete_dirichlet_condition} and~\eqref{eq:discrete_rigid_director} on $\gamma_0$.
\end{problem}

Note that frame indifference of the discrete model is retained naturally, because we simply restrict
the frame-indifferent functional $I$ to a subset $V_{p_1,h}^{\R^3} \times V_{p_2,h}^{\SOdrei}$ of its original domain of definition,
and this subset is closed under rigid body motions (Lemma~\ref{lem:equivariance_of_gfe}).

\begin{remark}
We have discretized the midsurface deformation $m$ using standard finite elements, and we have used the novel
geodesic finite elements only for the rotation field $\overline{R}$.  We can unify the two approaches when a more
abstract viewpoint is taken.  Indeed, revisiting the definitions of Chapter~\ref{sec:GFE} it is obvious that
geodesic finite elements may as well be defined for the target manifold $\R^3$ instead of $\SOdrei$; standard
Lagrangian finite elements are the result.  In this sense, we have used geodesic finite elements for both
the midsurface deformation and the microrotation field.

When the two orders $p_1$ and $p_2$ coincide $p=p_1 = p_2$, we can go one step further.  Note that the space $\text{SE(3)} \colonequals \R^3 \times \SOdrei$
is well known as the Special Euclidean group (the group of rigid body motions in $\R^3$).
We therefore introduce the GFE space $V_{h,p}^\text{SE(3)}$, and observe that it is isomorphic to $\mathbf{V}_h \colonequals V_{p,h}^{\R^3} \times V_{p,h}^\text{SO(3)}$.
We can therefore also interpret the discrete Cosserat shell problem as a minimization problem
in the single GFE space $V_{h,p}^\text{SE(3)}$.
\end{remark}

\subsection{The algebraic problem}

For the numerical minimization of the Cosserat shell energy we need an algebraic formulation.
For standard finite elements there is a bijective correspondence between finite element functions and
coefficient vectors, via the representation of the functions with respect to a basis.  For geodesic
finite elements, the situation is more involved.  Since GFE functions are continuous by definition,
we can always associate a coefficient vector $\overline{\overline{R}} \in \SOdrei^{N_2}$ to a function $\overline{R}_h \in V_{p_2,h}^{\SOdrei}$ by pointwise evaluation
at the $N_2$ Lagrange nodes.
To formalize this we introduce the evaluation operator
\begin{equation*}
 \mathcal{E}_{p_2} : V_{p_2,h}^{\SOdrei} \to \SOdrei^{N_2},
 \qquad
 \mathcal{E}_{p_2}(\overline{R}_h)_i = \overline{R}_h(n_i),
 \qquad
 i= 1,\dots,N_2,
\end{equation*}
where $n_i \in \omega$, $i=1,\dots,N_2$ are the Lagrange nodes of the $p_2$-order FE space on the grid.
However, for a given set of coefficients $\overline{\overline{R}} \in \SOdrei^{N_2}$ there may be more than
one GFE function that interpolates $\overline{\overline{R}}$.  This happens when the set of values violates the assumptions
of Theorems~\ref{thm:unique_minimizer_first_order} or~\ref{thm:well_posedness_vague} (depending
on the finite element approximation order $p_2$).

All geodesic finite element functions that do comply with the conditions of Theorems~\ref{thm:unique_minimizer_first_order} or~\ref{thm:well_posedness_vague}
element-wise can be identified with coefficient sets $\overline{\overline{R}} \in \text{SO(3)}^{N_2}$.
In most cases this situation can be achieved by making the grid fine enough.  This has been formalized
in~\cite[Thm.\,5.2]{sander:2013}, which we repeat here, adapted to the Cosserat shell problem.
\begin{theorem}
Let $\overline{R} : \omega \to \SOdrei$ be Lipschitz continuous
in the sense that there exists a constant $L$ such that
\begin{equation*}
 \dist(\overline{R}(x), \overline{R}(y)) \le L\norm{x-y}
\end{equation*}
for all $x,y \in \omega$.  Let $\mathcal{G}$ be a grid of $\omega$ and $h$
the length of the longest edge of $\mathcal{G}$.  Set $\overline{\overline{R}} = \mathcal{E}_{p_2}(\overline{R})$,
tacitly extending the definition of $\mathcal{E}_{p_2}$ to all continuous functions $\omega \to \SOdrei$.
For $h$ small enough, the inverse of $\mathcal{E}_{p_2}$ has only a single value in $V_{p_2,h}^{\SOdrei}$
for each $\widetilde{\overline{R}} \in \SOdrei^{N_2}$ in a neighborhood of $\overline{\overline{R}}$.
\end{theorem}
The restrictions posed by this theorem do not appear to pose any difficulties in practice.  We therefore
assume in the following that $\mathcal{E}_{p_2}$ is a (local) bijection.

Analogously to $\mathcal{E}_{p_2}$ we define the corresponding operator $\mathcal{E}_{p_1}$ doing point-wise evaluation
of functions in $V_{p_1,h}^{\R^3}$.  With these operators, it is straightforward to define the algebraic Cosserat shell energy
\begin{equation}
\label{eq:algebraic_shell_energy}
 \bar{I} \; : \R^{3N_1} \times \SOdrei^{N_2} \to \R,
 \qquad\qquad
 \bar{I}(\bar{m},\overline{\overline{R}}) \colonequals I(\mathcal{E}_{p_1}^{-1}(\bar{m}), \mathcal{E}_{p_2}^{-1}(\overline{\overline{R}})),
\end{equation}
where $I$ is the functional~\eqref{eq:finite_strain_energy}.
The algebraic Cosserat shell problem then is:
\begin{problem}[Algebraic Cosserat shell problem]
Find a pair $\bar{m} \in \R^{3N_1}$, $\overline{\overline{R}} \in \SOdrei^{N_2}$
that minimizes $\bar{I}$, subject to suitable boundary conditions.
\end{problem}
Implementation of Dirichlet boundary conditions
for the deformation $m_h$ is straightforward.  For the rotation field we again have the choice between leaving
the rotation free, or prescribing the transversal director vector $\overline{R}_3$ (rigid director prescription)
\begin{equation*}
 (\overline{\overline{R}}_i)_3
 =
 \frac{g'_d(n_i)}{\abs{g'_d(n_i)}}
\end{equation*}
for all Lagrange nodes $n_i$ on the Dirichlet boundary $\gamma_0$.

\begin{remark}
 If $N_1 = N_2 = N$ we can also interpret the functional \eqref{eq:algebraic_shell_energy} as being defined on the
 manifold $(\R^3 \times \SOdrei)^N$.
\end{remark}

It was mentioned in Chapter~\ref{sec:continuous_model} that the shell energy is frame-invariant in the
sense that
\begin{equation*}
 I(Qm, Q\overline{R}) = I(m,\overline{R}),
\end{equation*}
where $Q$ is any element of $\SOdrei$, acting on functions in $H^1(\omega,\R^3)$ and $W^{1,q}(\omega,\SOdrei)$
by pointwise multiplication.
By the equivariance property (Lemma~\ref{lem:equivariance_of_gfe}) of geodesic finite elements this frame
invariance does not get lost by discretization.
\begin{theorem}
 The algebraic energy functional $\bar{I}$ is frame-invariant in the sense that
\begin{equation*}
 \bar{I}(Q\bar{m}, Q\overline{\overline{R}}) = \bar{I}(\bar{m},\overline{\overline{R}}),
\end{equation*}
for all $Q \in \SOdrei$, which, by an abuse of notation, now acts on the components of $\bar{m}$ and $\overline{\overline{R}}$.
\end{theorem}
This sets the geodesic FE discretization apart from alternative approaches like~\cite{muench:2007, mueller:2009},
which do not have this property.

\section{Numerical minimization of the algebraic energy}
\label{sec:numerical_minimization}

All previous work on nonlinear shell elements has used the Newton method to solve the resulting nonlinear
systems of equations.  However, it is well known that this method converges only locally.  Therefore, a sequence
of loading steps is traditionally used to obtain a solution.  These loading steps have to be selected
carefully to make sure that the Newton solver converges at each loading step.  This selection of loading
steps can be tedious in practice.

For energy minimization problems there exist globalized versions of the Newton method, i.e., methods that
converge for any initial iterate, without using intermediate loading steps.  One such method is
the so-called trust-region method~\cite{conn_gould_toint:2000}, which replaces each Newton step with a
quadratic minimization problem on a convex set.  Under reasonable conditions, it degenerates to a standard
Newton method when close enough to a solution, and hence local quadratic convergence is recovered.

While the standard trust-region method works for energies defined on Euclidean spaces, a generalization to
energies on Riemannian manifolds has been introduced and investigated by \citet{absil_mahony_sepulchre:2008}.
This Riemannian trust-region method can be applied to the algebraic Cosserat energy~\eqref{eq:algebraic_shell_energy},
which is defined on the product manifold $\R^{3N_1} \times \SOdrei^{N_2}$. As an extension
of Newton's method, it shows locally quadratic behavior.  On the other hand, it can be shown to converge
globally without intermediate loading steps.

\subsection{Trust-region methods}

We briefly review the trust-region method for Euclidean spaces~\cite{conn_gould_toint:2000},
and then show how it can be generalized to functionals on a Riemannian manifold.
Consider a twice continuously differentiable functional
\begin{equation}
\label{eq:euclidean_functional}
 J : \R^N \to \R,
\end{equation}
supposed to be coercive and bounded from below.  Given any initial iterate $x^0 \in \R^N$, we want to find
a local minimizer of $J$.

\begin{figure}
 \begin{center}
  \begin{overpic}[width=0.8\textwidth]{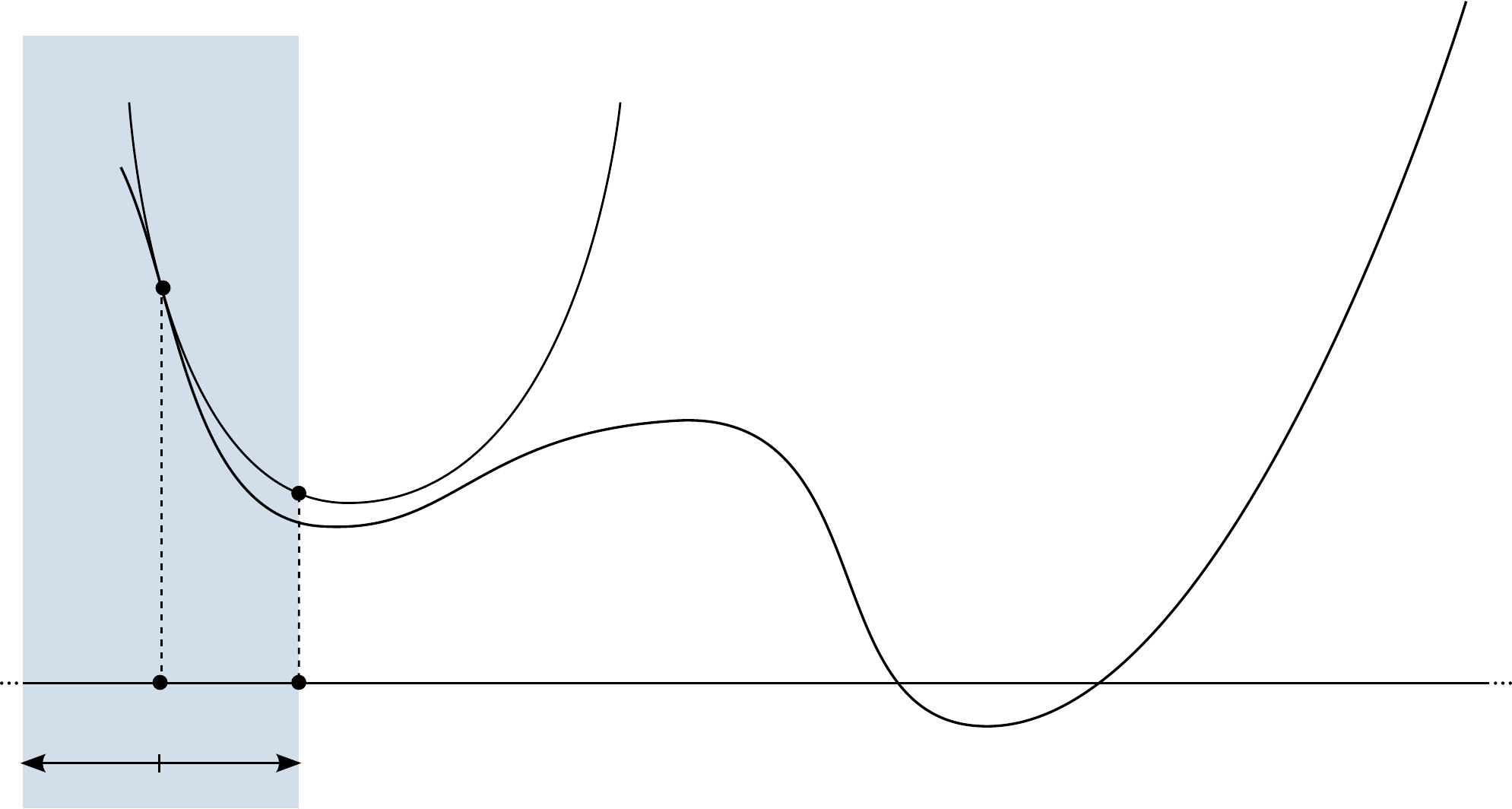}
  \put(9,5){$x^k$}
  \put(17,5){$x^{k+1}$}
  \put(41,40){$m_k$}
  \put(42,27){$J$}
  \put(5,1){$\rho_k$}
  \put(14,1){$\rho_k$}
  \put(95,4){$\R^N$}
  \end{overpic}
  \caption{One step of the trust-region method.  The new iterate $x^{k+1}$ is the minimizer of the quadratic
    model $m_k$ restricted to the ball $B_{x^k}(\rho_k)$ (shaded region), unless the energy decrease predicted by the model
    deviates too much from the true energy decrease $J(x^k) - J(x^{k+1})$.
   \label{fig:trust_region_schematic}
    }
 \end{center}
\end{figure}
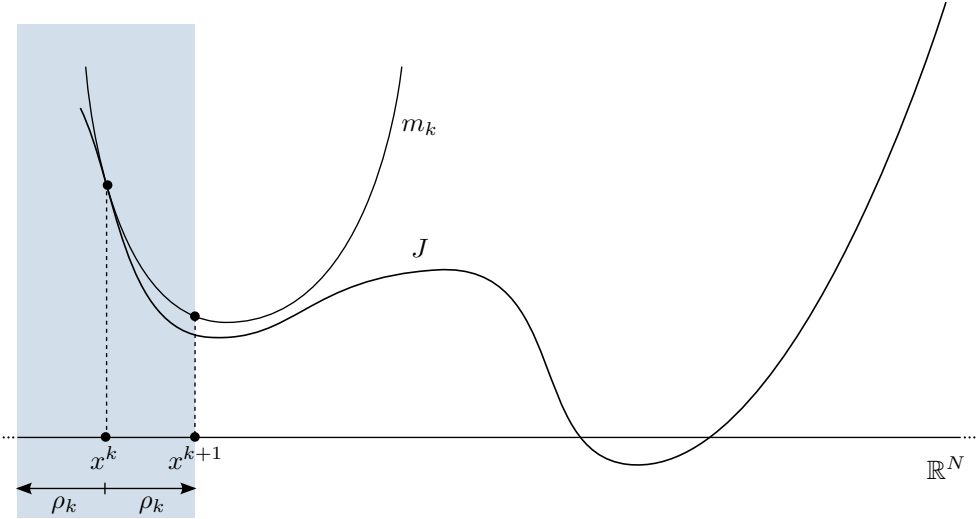

The Newton method does this in the following way.  Let $x^k \in \R^N$ be any iterate.  Approximate $J$
around $x^k$ by the quadratic Taylor expansion
\begin{align*}
 m_k & \; : \R^N \to \R, \\
 m_k(s) & = J(x^k) + \partial J(x^k)s + \frac{1}{2} s^T \partial^2 J(x^k) s,
\end{align*}
which in this context is called a quadratic model of $J$ around $x^k$.  The variable $s$ is to be
interpreted as a correction $s = x - x^k$.  Then, compute a stationary point $s^k$ of $m_k$, and use it
as the correction to the next iterate
\begin{equation*}
 x^{k+1} \colonequals x^k + s^k.
\end{equation*}
Computing the stationary point $s^k$ is done by the well-known Newton update formula
\begin{equation}
\label{eq:newton_step}
 s^k
 =
 x^{k+1}-x^k
 =
 - \partial^2 J(x^k)^{-1} \partial J(x^k).
\end{equation}
Observe that if the Hessian $\partial^2 J(x^k)$ is positive definite at all iterates, then the algorithm
produces a sequence of iterates with decreasing energy, i.e., $J(x^{k+1}) \le J(x^k)$ for all $k \in \mathbb{N}$.
However, iterates with indefinite $\partial^2 J(x^k)$ may lead to energy increase.

To enforce global convergence of this, the trust-region method first replaces the search
for a stationary point of $m_k$ by a minimization problem for a minimizer $s^k$ of $m_k$.  As a consequence,
iterates of the trust-region method are energy decreasing in all cases.  Secondly, it notes that the quadratic model $m_k$
is a good approximation of $J$ only in a neighborhood of $x^k$.  This observation is made explicit
by restricting the minimization problem for $m_k$ to a ball of radius $\rho_k$ around $x^k$,
the name-giving trust region (Figure~\ref{fig:trust_region_schematic}).  In other words, the Newton step~\eqref{eq:newton_step} is replaced by
\begin{equation}
\label{eq:trust_region_step}
 s^k = \argmin_{s \in \R^N} m_k(s),
 \qquad
 \norm{s^k} \le \rho_k.
\end{equation}
Since we now look for a minimizer on a compact set only,
Problem~\eqref{eq:trust_region_step} is well-defined even if $\partial^2 J$ is not positive definite.

Unlike the original Newton method, the trust-region method is monotone in the sense
that $J(x^{k+1}) \le J(x^k)$ for all $k \in \mathbb{N}$.  A more quantitative monitoring of the
energy decrease allows to control the trust-region radius, i.e., the trust in the quality of the
quadratic approximation.
The quality of the correction step $s^k$ is estimated by comparing the functional decrease to the model decrease.
If the quotient
\begin{equation}
\label{eq:tr_descent_control}
  \kappa_k = \frac{J(x^k) - J(x^k + s^k)}{m_k(0) - m_k(s^k)}
\end{equation}
is smaller than
a fixed value $\eta_1$, then the step is rejected, and $s^k$ is recomputed for a
smaller trust-region radius.  Otherwise the step is accepted.
If $\kappa_k$ is larger than a second value $\eta_2$,
the trust-region radius is enlarged for the next step.  Common values are $\eta_1 = 0.01$ and
$\eta_2 = 0.9$~\cite{conn_gould_toint:2000}.

For the trust-region algorithm, the following convergence properties can be shown.
\begin{theorem}[{ \cite[Thms.\,6.4.6 and 6.5.5]{conn_gould_toint:2000}}]
\label{thm:tr_convergence}
 Suppose that $J$ is twice continuously differentiable, bounded from below, and such that its Hessian
 remains bounded for all $x \in \R^N$.
 \begin{enumerate}
 \item For all initial iterates we get
 \begin{equation*}
  \lim_{k \to \infty} \norm{\partial J(x^k)} = 0.
 \end{equation*}
 \item Suppose that $\{x^{k_i}\}$ is a subsequence of the iterates
 converging to the first-order critical point $x_*$.  Suppose furthermore that
 $s^k \neq 0$ for all $k$ sufficiently large.  Finally suppose that $\partial^2 J(x_*)$ is positive definite.
 Then the complete sequence of iterates $\{x^k\}$ converges to $x_*$, eventually the step quality $\kappa_k$
 remains above $\eta_2$,  and the trust-region radius $\rho_k$ is bounded away from zero.
 \end{enumerate}
\end{theorem}

In particular, since $\kappa_k > \eta_2$ for all $k$ large enough, the trust-region radius grows near local
minimizers, the method eventually degenerates to a pure
Newton method, and we get locally quadratic convergence.

Various algorithms for solving the constrained quadratic minimization problems \eqref{eq:trust_region_step} been proposed
in the literature.  The monograph~\cite{conn_gould_toint:2000} gives a good overview.

\subsection{Riemannian trust-region methods}
\label{sec:riemannian_trust_region}

\begin{figure}
 \begin{center}
  \begin{overpic}[width=0.47\textwidth]{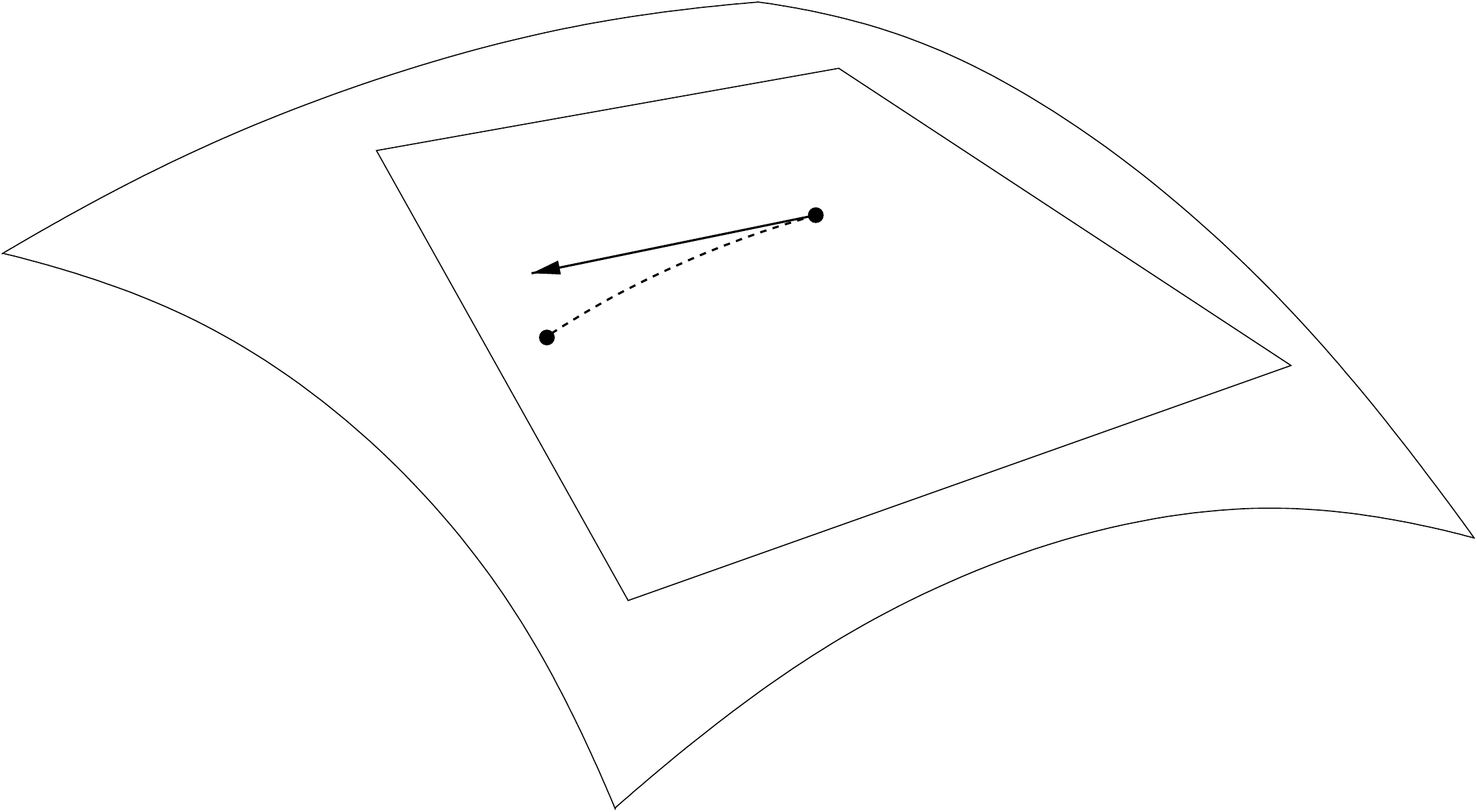}
  \put(56,37){$x^k$}
  \put(45,40){$s^k$}
  \put(38,27){$x^{k+1}$}
  \put(10,36){$M$}
  \put(68,29){$T_{x^k}M$}
  \end{overpic}
 \end{center}
  \caption{In the Riemannian trust-region method, the energy functional defined on $M$ is
    lifted onto the tangent space at $x^k$ using the exponential map.  Then, a linear
    correction step is computed on $T_{x^k} M$, and applied to $x^k$ using the exponential
    map $x^{k+1} = \exp_{x^k} s^k$.}
  \label{fig:retraction}
\end{figure}
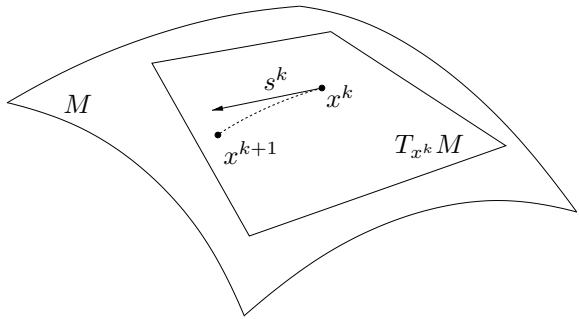

The algebraic energy functional $\bar{I}$ defined in \eqref{eq:algebraic_shell_energy} is not a functional of the type~\eqref{eq:euclidean_functional}.
Rather, its domain of definition is the nonlinear manifold $\R^{3N_1} \times \SOdrei^{N_2}$.
The trust-region method has been generalized to such energies by~\citet{absil_mahony_sepulchre:2008}.
Let $M$ be a Riemannian manifold with metric $g$, and $J : M \to \R$ twice
differentiable and bounded from below (in our case: $M =\R^{3N_1} \times \SOdrei^{N_2}$). The basic idea of such a Riemannian trust-region
algorithm is that in a neighborhood of a point $x \in M$ the functional $J$ can be lifted
onto the tangent space $T_xM$.  There, a vector space trust-region subproblem can be solved and
the result transported back onto $M$ (Figure~\ref{fig:retraction}).

More formally, let again $k \in \mathbb{N}$ be an iteration number and $x^k \in M$ the current iterate.
We obtain the lifted functional by setting
\begin{align*}
\hat{J}_k  \; : \; T_{x^k} M \to \R,
\qquad\qquad
\hat{J}_k (s) = J(\exp_{x^k}s).
\end{align*}
Let $\rho_k > 0$ be the current trust-region radius.
The Riemannian metric $g$ turns $T_{x^k} M$ into a Banach space with the norm
$\norm{\cdot}_{x^k} = \sqrt{g_{x^k}(\cdot, \cdot)}$.  There, the trust-region
 subproblem reads
\begin{equation}
\label{eq:rtr_subproblem}
s_k = \argmin_{s \in T_{x^k}M} m_k(s),
\qquad
\norm{s}_{x^k}  \le \rho_k,
\end{equation}
with the quadratic, but not necessarily convex model
\begin{equation}
\label{eq:rtr_quadratic_functional}
m_k(s)
=
\hat{J}_k(0) + g_{x^k}(\nabla \hat{J}_k(0),s)
 + \frac 12 g_{x^k} ( \hessian \hat{J}_k(0)s,s).
\end{equation}
Here
$\nabla \hat{J}_k$ is the Riemannian gradient
and $\hessian \hat{J}_k$ the Riemannian Hessian of $\hat{J}_k$ (see~\cite{absil_mahony_sepulchre:2008} for definitions),
and both are evaluated at $0 \in T_{x^k}M$.
Note that \eqref{eq:rtr_quadratic_functional} is independent of a specific coordinate system on $T_{x^k}M$.
As a minimization problem of a continuous function on a compact set,
\eqref{eq:rtr_subproblem} has at least one solution $s^k$, which
generates the new iterate by
\begin{equation*}
 x^{k+1} = \exp_{x^k} s^k.
\end{equation*}

As in trust-region methods in linear spaces, the quality of a correction step $s^k$
is estimated by comparing the functional decrease and the model decrease.
The quotient~\eqref{eq:tr_descent_control} now takes the form
\begin{equation*}
  \kappa_k = \frac{J(x^k) - J(\exp_{x^k} s^k)}{m_k(0) - m_k(s^k)}.
\end{equation*}

For this method, \citeauthor{absil_mahony_sepulchre:2008} proved global convergence to first-order
stationary points, and, depending on the exactness of the inner solver,
locally superlinear or even locally quadratic convergence~\cite{absil_mahony_sepulchre:2008}.
For our numerical
results we use the monotone multigrid method~\cite{kornhuber:1997} together with a $\infty$-norm trust-region.
Details can be found in~\cite{sander:2012}.

\subsection{Computing the algebraic tangent problem numerically}
\label{sec:computing_tangent_problem}

Solving the constrained quadratic problems~\eqref{eq:rtr_subproblem} numerically involves the algebraic
Riemannian gradient $\nabla \bar{I}$ and Hessian $\hessian \bar{I}$ of the functional $\bar{I}$.  While those could in principle be evaluated
analytically, such an approach is involved and error prone
(Consider the derivative formulas for the gradient in~\cite[Chap.\,5]{sander:2012}).  It is much more convenient to use automatic
differentiation (AD) to compute the derivatives.
AD is a technique to algorithmically compute first and higher derivatives of functions given in form of
computer programs~\cite{griewank_walther:2008}.  This includes computer programs involving iterative solvers like the Newton method
used to evaluate GFE functions (Section~\ref{sec:values_of_gfe_functions}).  Many good implementations of AD
are available as external libraries.  For this article we have used the open-source ADOL-C software~\cite{walther_griewank:2012}.

For the rest of this paper we assume that the deformation $m$ and the microrotation $\overline{R}$ have been discretized with
finite elements of equal approximation order.  Then there is an an equal number of Lagrange nodes $N = N_1 = N_2$ for both of them,
and we can consider the algebraic energy $\bar{I}$ as being defined on the manifold $M = (\R^3 \times \SOdrei)^N$.

Unfortunately, current AD tools do not directly support derivatives of energies defined on manifolds.
We therefore use the following trick.  Interpret elements $R$ of $\SOdrei$ as unit vectors $q$ in $\R^4$
using quaternion coordinates (see Section~\ref{sec:quaternion_coordinates}).
The algebraic energy functional $\bar{I}$ can then be interpreted as being defined on $(\R^3 \times S^3)^N \subset \R^{7N}$.
To extend $\bar{I}$ to a neighborhood of $(\R^3 \times S^3)^N$ in $\R^{7N}$ we first introduce $\bar{q} \in \R^{4N}$,
a vector of quaternions.  Componentwise normalization leads to a vector of unit quaternions, which we denote by
$\bar{q} / \abs{\bar{q}} \in (S^3)^N$ in an abuse of notation.  Using the map $F$ defined in~\eqref{eq:SO3_double_cover} we
can construct $F(\bar{q}/\abs{\bar{q}}) \in \SOdrei^N$ (the application of $F$ again component-wise).  Then we set
\begin{equation*}
 \tilde{I} (\bar{m}, \bar{q})
 \colonequals
 \bar{I} \big (\bar{m}, F(\bar{q} / \abs{\bar{q}}) \big),
\end{equation*}
which is a smooth functional on an open subset of the Euclidean space $\R^{7N}$.
Given a computer implementation of $\tilde{I}$,
an AD system like \textsc{ADOL-C} can then compute the Euclidean gradient $\partial \tilde{I} \in \R^{7N}$ and
Hessian $\partial^2 \tilde{I} \in \R^{7N \times 7N}$ automatically.

To obtain the Riemannian gradient $\nabla \bar{I}$ and Hessian $\hessian \bar{I}$ we need additional manipulations.
For the gradient we use the following well-known result (see, e.g., \cite{absil_mahony_sepulchre:2008}, Sec.\,3.6.1).
\begin{lemma}
\label{lem:embedded_gradient}
Let $M$ be a smooth Riemannian
manifold isometrically embedded in a Euclidean space $\R^l$.
For each $x \in M$ let $P_x : T_x\R^l \to T_xM$ be the orthogonal projection onto
the tangent space at $x$.
Let $f : M \to \R$ be continuously differentiable and $\tilde{f}$ a smooth
extension of $f$ to a neighborhood of $M$ in $\R^l$.  Then
\begin{equation}
 \label{eq:projected_gradient}
 \nabla f = P_x \partial \tilde{f},
\end{equation}
where $\nabla$ is the gradient operator on $M$, and $\partial$ is the gradient in $\R^l$.
\end{lemma}
Since $\bar{I}$ is defined on the $N$-fold product of $\R^3 \times \SOdrei$ we obtain the Riemannian
gradient $\nabla \bar{I}$ by applying Lemma~\ref{lem:embedded_gradient} to each factor.
Hence, the Riemannian gradient is given by componentwise projection
\begin{equation*}
( \nabla \bar{I})_i
 =
 P_x (\partial \tilde{I})_i,
 \qquad
 i = 1,\dots, N,
\end{equation*}
where $P_x$ is the orthogonal projector from $v \in \R^7$ to $\R^3 \times T_x S^3$.  This projector can be constructed
from the corresponding projector for $\R^3$ (which is the identity), and the corresponding projector for $S^3$
\begin{equation*}
 P^{S^3}_x = I - x x^T.
\end{equation*}

A similar formula for the Riemannian Hessian is given in the following lemma.  As we now consider
second derivatives, the curvature of $\SOdrei$ comes into play.
\begin{lemma}[\citet{absil_mahony_trumpf:2013}]
\label{lem:embedded_hessian}
With the same notation as in Lemma~\ref{lem:embedded_gradient}, we have
\begin{equation*}
 \hessian f(x)[z]
 =
 P_x \partial^2 \tilde{f}(x)z + \mathfrak{A}_x(z,P^\perp_x \partial \tilde{f}),
\end{equation*}
where $\mathfrak{A}_x(z,v)$ is the Weingarten map of $M$, and $P_x^\perp$
is the orthogonal projector onto the normal space of $M$ at $x$.
\end{lemma}
The Weingarten map for the unit sphere in $\R^4$ is~\cite{absil_mahony_trumpf:2013}
\begin{equation*}
 \mathfrak{A}_x(z,v)
 \colonequals
 - (x^Tv)z,
\end{equation*}
and the orthogonal projector onto the normal space at $x \in S^3$ is
\begin{equation*}
 P_x^\perp
 =
 I - P_x
 =
 xx^T.
\end{equation*}

\bigskip

Written in canonical coordinates of $\R^{7N}$, the matrix $\hessian \tilde{I}$ is a sparse symmetric
$7N \times 7N$-matrix, consisting of dense $7 \times 7$ blocks.  Using this representation for numerical
computations is undesirable for two reasons.   First of all, it is rank deficient, because the extended
functional $\tilde{I}$ is constant along each normal vector of $S^3$.  Secondly, it is bigger than
necessary: Since $\SOdrei$ (or the set of unit quaternions for that matter) is only three-dimensional, the entire Riemannian
Hessian should fit into a $6N \times 6N$ matrix.
To construct such a representation for the Riemannian Hessian at a point $(\bar{m}, \bar{q}) \in \R^{7N}$
we pick a basis for the tangent space of $(\R^3 \times S^3)^N$ at $(\bar{m}, \bar{q})$, and write
$\hessian \tilde{I}$ in that basis.  Luckily, such a basis is easily available.
For the components in $\R^3$, the canonical basis can be used.  For any point $q \in S^3$, an orthonormal
basis of $T_q S^3$ is given by
\begin{equation*}
D_{q,1} 
=
\begin{pmatrix}
  q_3 \\  q_2  \\ -q_1 \\ -q_0
\end{pmatrix},
\qquad
D_{q,2} 
=
\begin{pmatrix}
 -q_2 \\  q_3 \\  q_0 \\ -q_1
\end{pmatrix},
\qquad
D_{q,3} 
=
\begin{pmatrix}
  q_1 \\ -q_0 \\  q_3 \\ -q_2
\end{pmatrix},
\end{equation*}
and this basis depends smoothly on $q$.  We combine the vectors to a $7 \times 6$-matrix
\begin{equation}
\label{eq:tangent_space_basis}
D_q
=
\begin{pmatrix}
1 & & & & & \\
  & 1 & & & & \\
  & & 1 & & & \\
  & & & q_3 & -q_2 & q_1 \\
  & & & q_2 & q_3 & -q_0 \\
  & & & -q_1 & q_0 & q_3 \\
  & & & -q_0 & -q_1 & -q_2
\end{pmatrix},
\end{equation}
whose columns form an orthonormal basis of $\R^3 \times S^3$.

We denote by $D$ the block-diagonal $7N \times 6N$-matrix where the $i$-th block is $D_{q_i}$ as given by~\eqref{eq:tangent_space_basis}.
Then, in these new coordinates, the Riemannian Hessian has the algebraic form
\begin{equation}
\label{eq:riemannian_hessian}
 \hessian \tilde{I}
 =
 D^T \partial^2 \tilde{I} D - D^T(x^T P^\perp_x \partial \tilde{I}) D
 \quad \in
 \R^{6N \times 6N}.
\end{equation}
This matrix has no degenerate directions caused by the embedding of the configuration space into $\R^{7N}$.
Indeed, it is again completely intrinsic.  In each iteration of the trust-region solver, this is the matrix
used to define the quadratic model.

Finally, we point out one lucky coincidence that helps to increase efficiency. AD systems such as ADOL-C
are able to compute the product $(\partial^2 \tilde{I}) D$ directly.  This is noticeably cheaper
than using AD to compute $\partial^2 \tilde{I}$ and later multiplying by $D$, because $(\partial^2 \tilde{I}) D$ has fewer entries than $\partial^2 \tilde{I}$
($7N \times 6N$ compared to $7N \times 7N$).
We noted a decrease of about 10\% of the time needed to assemble the Riemannian Hessian~\eqref{eq:riemannian_hessian}.

\section{Numerical tests}
\label{sec:numerical_tests}

We now present several numerical tests.  These demonstrate the capabilities
of both our Cosserat shell model and of our discretization.  We reproduce quantitative results from the literature~(Section~\ref{sec:wriggers-l-shape}),
and show how the model and discretization can handle large rotations with ease (Section~\ref{sec:twisted_strip}).  In Section~\ref{sec:wrinkling}
we simulate the wrinkling of a polyimide sheet, and find very good quantitative correspondence with experimental data.
All examples in this chapter were
programmed using the \textsc{Dune} libraries (\cite{bastian_et_al:2008}, \url{www.dune-project.org}).
For all examples we used second-order finite elements for both the deformation $m$ and the microrotation $\overline{R}$.
No locking effects could be observed for this discretization.

\subsection{Deformation of an \texorpdfstring{$L$}{L}-shape}
\label{sec:wriggers-l-shape}

\begin{figure}
 \begin{center}
  \begin{overpic}[width=0.9\textwidth]{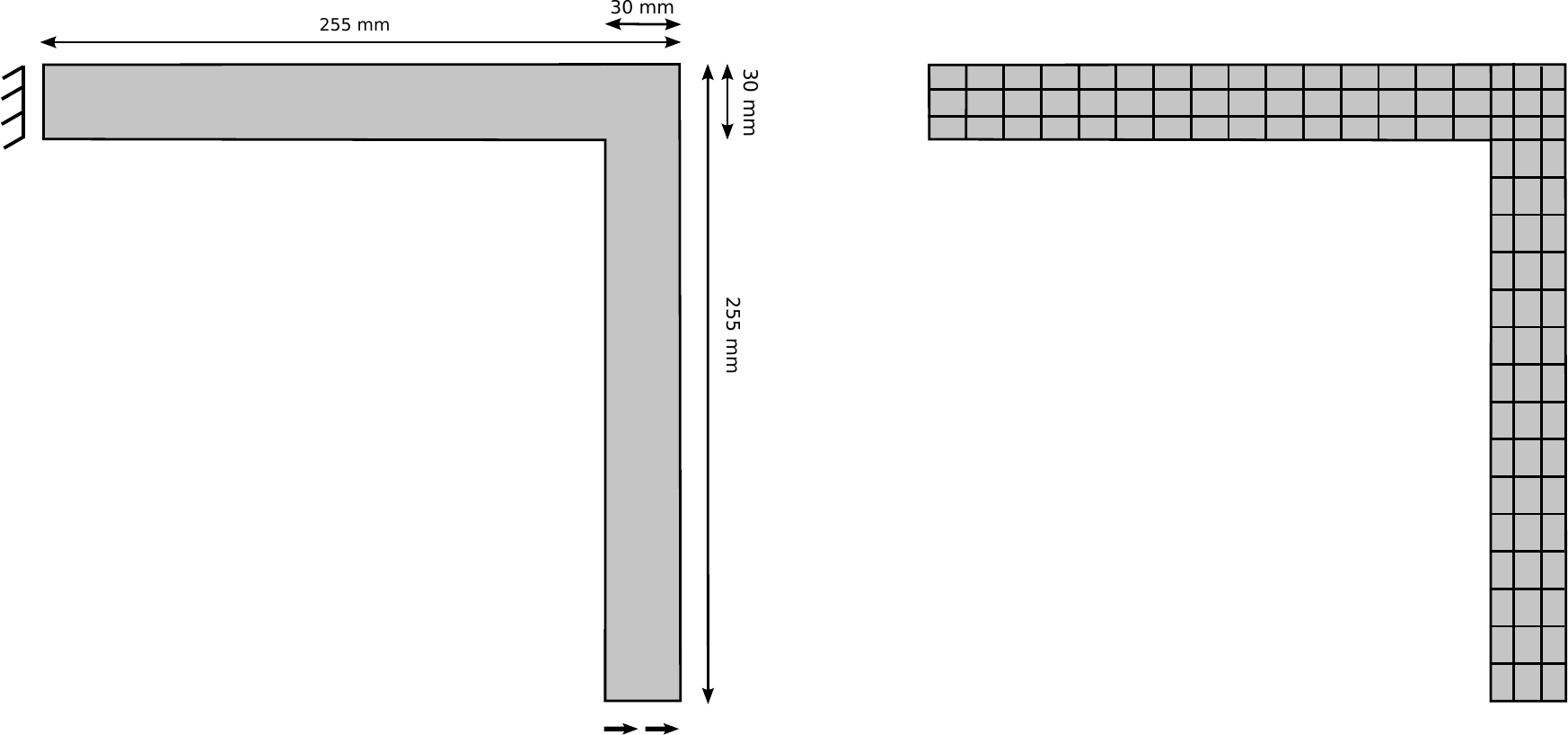}
   \put(3.5,40){$\gamma_0$}
   \put(40,3.5){$\gamma_s$}
  \end{overpic}
 \end{center}
 \caption{Left: $L$-shape structure with boundary conditions; Right: the grid,
    which is the one also used in~\cite{wriggers_gruttmann:1993}}
 \label{fig:L-shaped_domain}
\end{figure}

We begin by comparing our approach to a benchmark problem taken from the literature.  The following setup
is used by \citet{wriggers_gruttmann:1993}, who compare their discrete model with the ones from
\cite{argyris_et_al:1979,simo_vu-quoc:1986,simo_fox_rifai:1990} for the same problem.
Our aim here is two-fold: We want to show that our discrete model can reproduce quantitative results from the literature.
Also, we want to highlight the speed and stability of our solver.

Let $\omega$ be the $L$-shaped domain depicted in Figure~\ref{fig:L-shaped_domain}.
Sizes of the shape are given in the figure, and we set the plate thickness to $0.6$\,mm.
We model the material with the finite-strain hyperelastic material of Section~\ref{sec:finite_strain_shell_model}.
The material parameters are given in Table~\ref{tab:material_parameters_L-shape}.
The Lam\'e constants $\mu,\lambda$ correspond to the values $E = 71\,240\,\text{N/mm}^2$, $\nu=0.31$ given
in~\cite{wriggers_gruttmann:1993}.  As argued in Section~\ref{sec:infinitesimal_shell_model},
the coupling modulus $\mu_c$ is set to $\mu_c = 0$\,N/mm.
We set the curvature exponent $q$ appearing in the curvature energy term $W_\text{curv}$ to $q=2$,
and the internal length $L_c$ to $0.6\,\mu$m, following the suggestions of Section~\ref{sec:infinitesimal_shell_model}.

The boundary conditions are depicted on the left of Figure~\ref{fig:L-shaped_domain}.
The structure is clamped on the left vertical end $\gamma_0$.  By this we mean that on $\gamma_0$ we set
$m(x,y) = (x,y,0)$, and the rigid director description $\overline{R}_3 = (0,0,1)^T$ for the
microrotations $\overline{R}$. On the lower horizontal end $\gamma_s$ we prescribe a uniform surface load%
\footnote{Here we deliberately differ from~\cite{wriggers_gruttmann:1993}, where a point load is used.}
$P$ in the direction of the first unit basis vector.
Zero Neumann boundary conditions are set everywhere else for displacements and rotations.
We discretize the domain using 99 quadrilateral elements as depicted on the right of Figure~\ref{fig:L-shaped_domain}.
The equations are discretized using second-order (i.e., nine-node) geodesic finite elements.
\begin{table}
\begin{center}
 \begin{tabular}{cccccc}
  $h\,[\text{mm}]$ 
  & $\mu\,[\text{N}/\text{mm}^2]$ & $\lambda\,[\text{N}/\text{mm}^2]$ 
  & $\mu_c\,[\text{N}/\text{mm}^2]$ 
  & $L_c\,[\text{mm}]$   
  & $q\,[1]$     
  \\
  \hline
  0.6 & $2.7191\cdot 10^{4}$ & $4.4364 \cdot 10^{4} $ & 0 & $0.6\cdot 10^{-3}$ & 2
 \end{tabular}
\end{center}
\caption{Material parameters for the $L$-shape example}
\label{tab:material_parameters_L-shape}
\end{table}

\begin{figure}
 \begin{center}
  \includegraphics[width=0.9\textwidth]{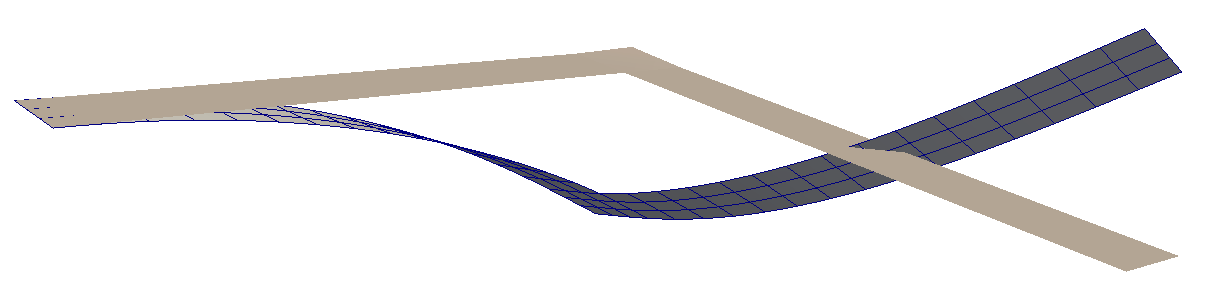}
  \includegraphics[width=0.7\textwidth]{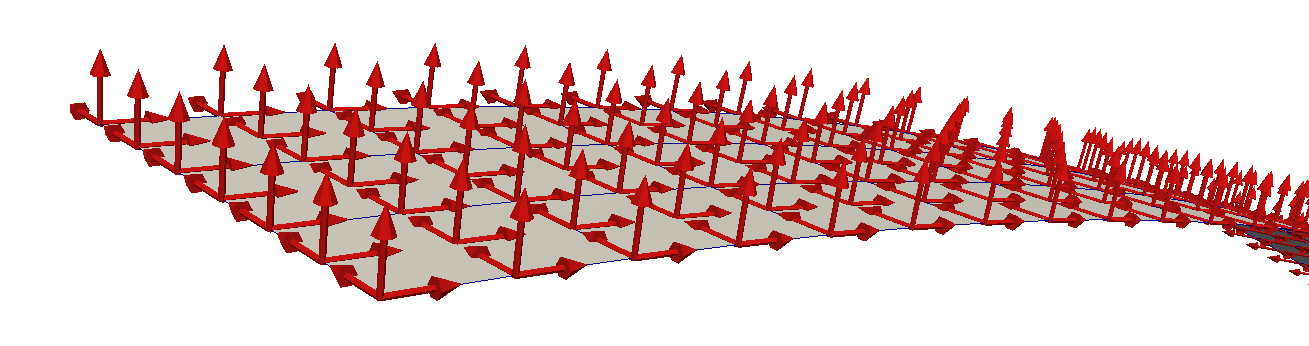}
 \end{center}
 \caption{Example deformation of the $L$-shape structure for $P=1.62$\,N.  Upper picture: initial configuration and
 configuration under load. Lower picture: closeup of the clamped part of the structure,
 with the directors shown as red arrows}
 \label{fig:L-shape_example_deformation}
\end{figure}

\bigskip

The first aim of this experiment is to study the buckling behavior of the structure for different values of $P$.
When the structure is loaded, it deforms in-plane as long as the load $P$ stays below a critical value $P_s$.
For loads beyond this value, the structure starts to buckle laterally.
An example deformation using $P = 1.62$\,N is shown in Figure~\ref{fig:L-shape_example_deformation}.

Since the in-plane deformation remains a stationary point of the energy even for loads larger than $P_s$,
a perturbation needs to be applied to trigger the buckling.  We do this by starting the trust-region method
at the asymmetric initial iterate
\begin{equation}
 \label{eq:l-shape-initial-iterate}
 m(x,y) = \Bigg(x,y,z=
 \begin{cases}
  0 & \text{if $x < 225$ or $y < -15$} \\
  10^{-3}(x-225)(y+15) & \text{else}
 \end{cases}
 \Bigg)
 \qquad
 \overline{R} = \identity.
\end{equation}
This adds a little kink in the corner of the domain, which is enough to trigger the buckling.

\begin{figure}
 \begin{center}
  \includegraphics[width=0.42\textwidth]{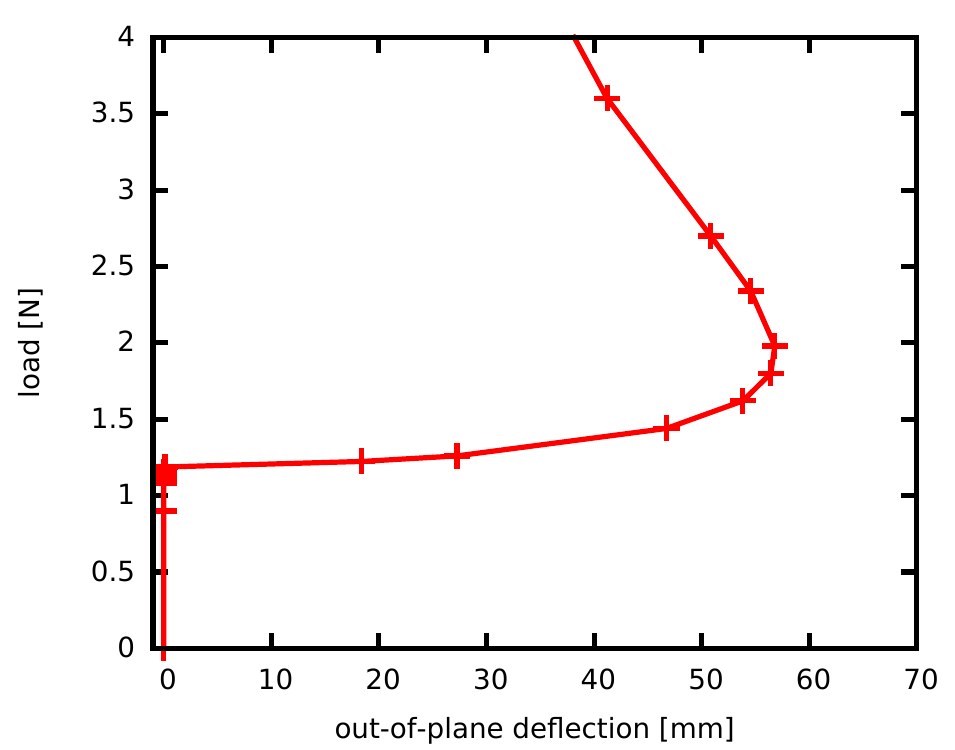}
  \quad
  \includegraphics[width=0.42\textwidth]{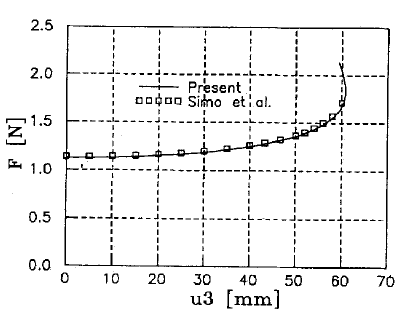}
 \end{center}
 \caption{Out-of-plane deflection as a function of the load.  Left: own simulation.  Right: corresponding plot taken from~\cite{wriggers_gruttmann:1993}}
 \label{fig:out_of_plane_deflection}
\end{figure}

A plot showing the lateral average displacement of $\gamma_s$ is shown in Figure~\ref{fig:out_of_plane_deflection}.
For comparison we have also given the corresponding plot from~\cite{wriggers_gruttmann:1993}.
It can be seen that the critical value we obtain is between 1.188\,N and 1.224\,N.
This is in good agreement to the other values from the literature\cite{wriggers_gruttmann:1993,simo_vu-quoc:1986,argyris_et_al:1979, simo_fox_rifai:1990},
which we print in Table~\ref{tab:l-shape-stability}.

\bigskip

\begin{table}
\begin{center}
 \begin{tabular}{cccccc}
 Reference & Type & Type of elements & Number of elements & $P_s$\,[N] \\
 \hline
 \cite{argyris_et_al:1979}       & beam  & ---           & 20 & 1.088 \\
 \cite{simo_vu-quoc:1986}        & beam  & ---           & 20 & 1.090 \\
 \cite{argyris_et_al:1979}       & shell & triangle      & 86 & 1.145 \\
 \cite{simo_fox_rifai:1990}      & shell & quadrilateral & 68 & 1.137 \\
 \cite{simo_fox_rifai:1990}      & shell & quadrilateral & converged solution & 1.128\\
 \cite{wriggers_gruttmann:1993}  & shell & nine-node     & 17 & 1.113 \\
 \cite{wriggers_gruttmann:1993}  & shell & nine-node     & 99 & 1.123 \\
 here                            & shell & nine-node     & 99 & 1.188--1.224
 \end{tabular}
\end{center}
\caption{Literature results for the critical load}
\label{tab:l-shape-stability}
\end{table}

In a second step we want to highlight a few properties of the solver.  For this we use the configuration described above
with the surface load $P= 1.62$\,N at $\gamma_s$ shown in Figure~\ref{fig:L-shape_example_deformation}.
We solve the problem in a single loading step, using the trust-region method described in Section~\ref{sec:riemannian_trust_region}.
For the quadratic minimization problems we use a monotone multigrid method as described in~\cite{sander:2012}.
The $\infty$-norm is used to define the trust region.  We scale the rotation part of the norm by a factor of $10^{-3}$,
so that corrections to the deformation (with numerical values in the two-digit range) are treated equally to corrections
to the rotations (which cannot get larger than $\pi$).

\begin{figure}
 \begin{center}
  \includegraphics[width=0.32\textwidth]{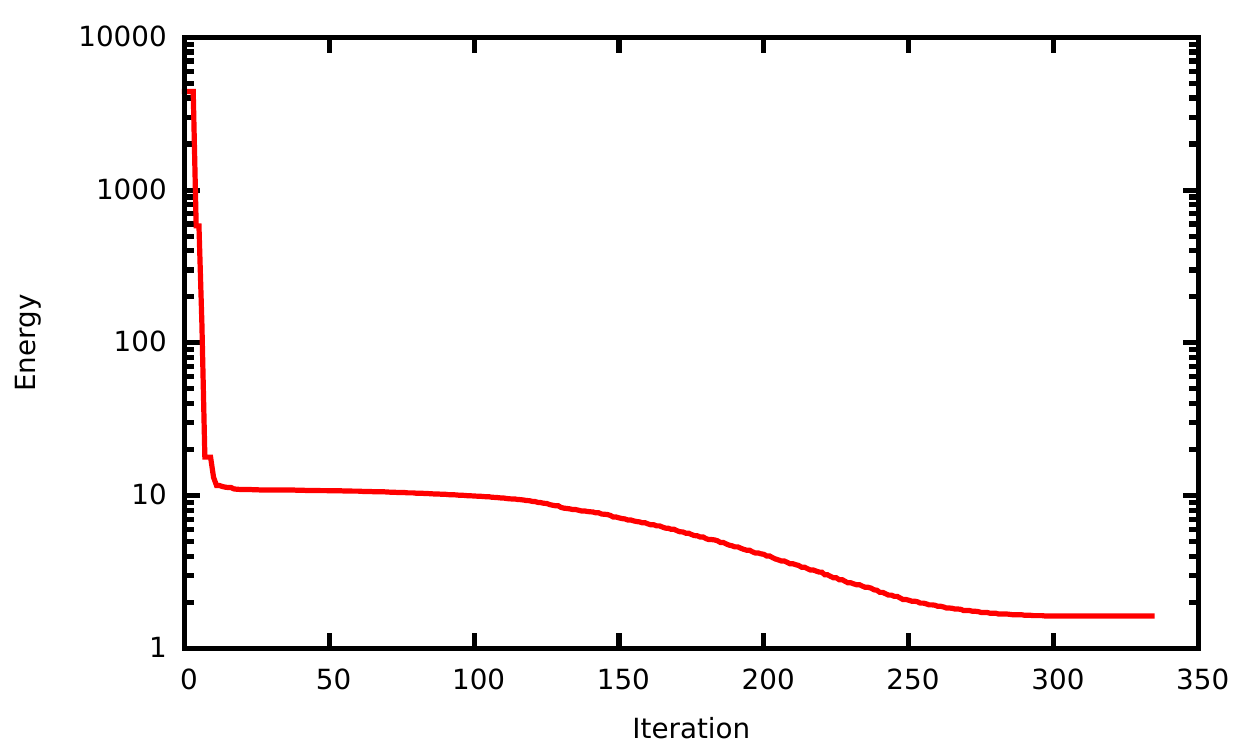}
  \includegraphics[width=0.32\textwidth]{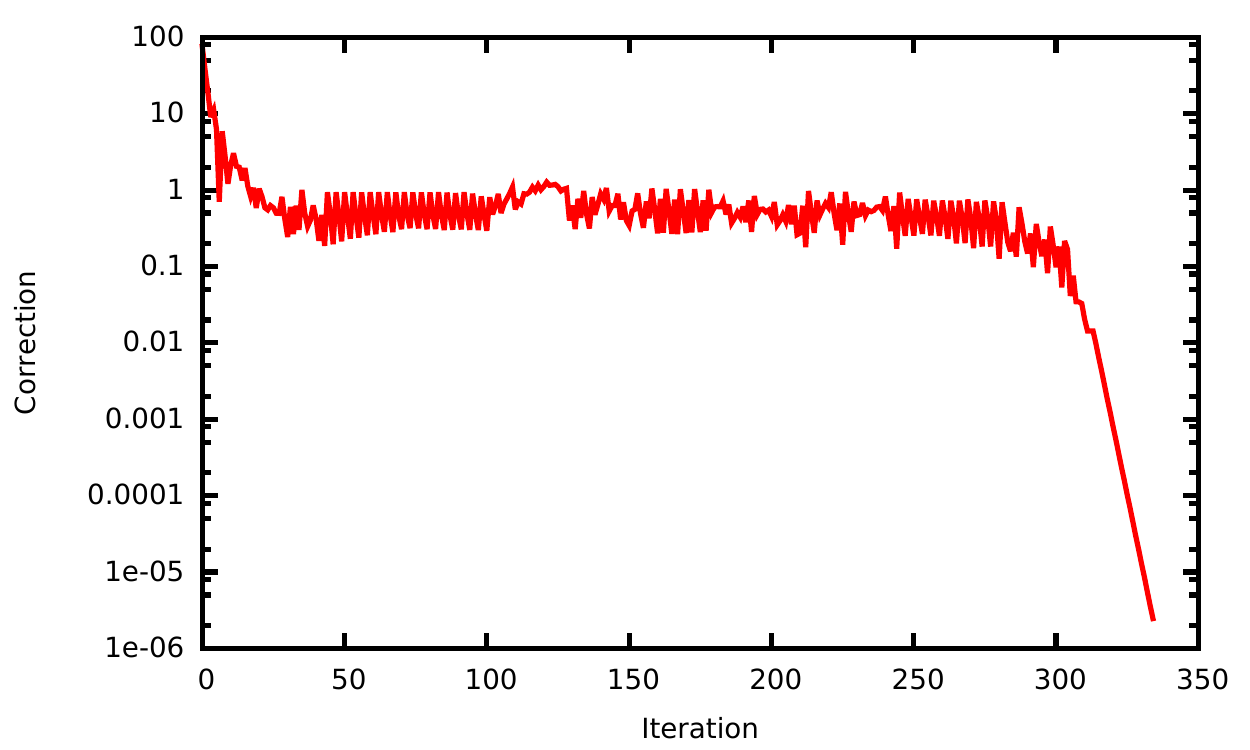}
  \includegraphics[width=0.32\textwidth]{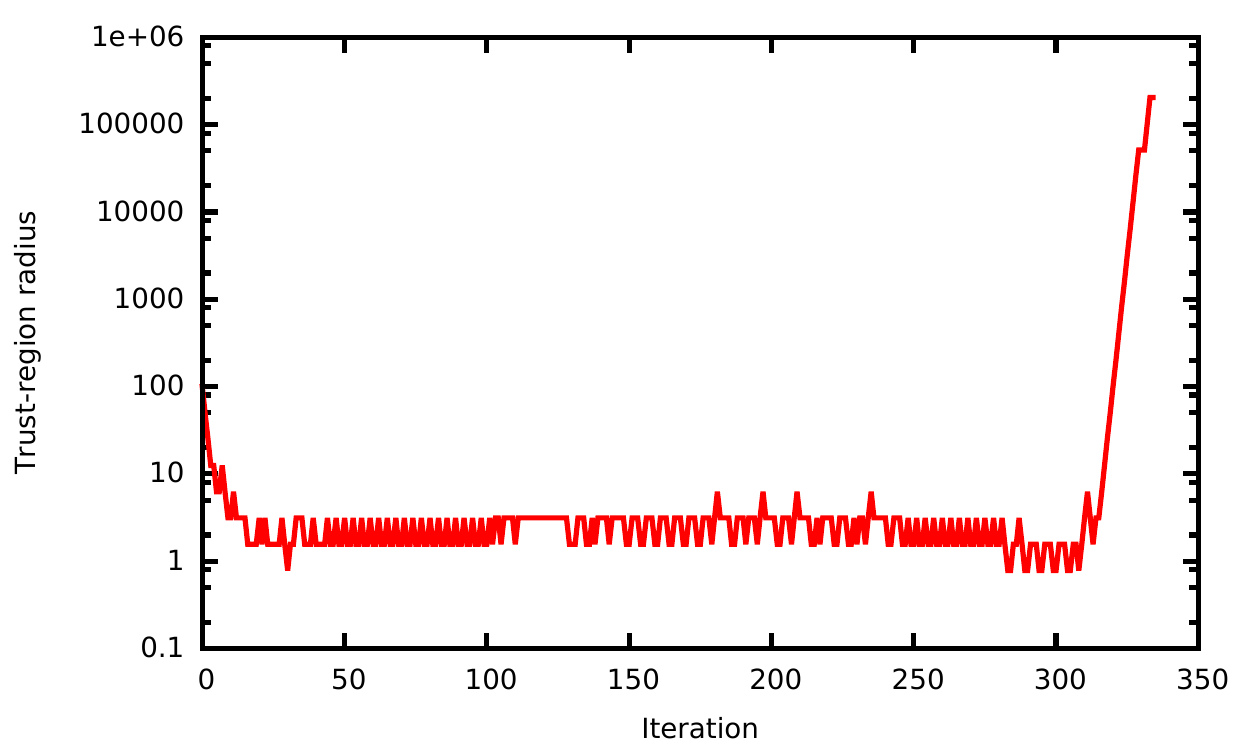}
 \end{center}
 \caption{Behavior of the Riemannian trust-region solver for the configuration shown in
   Figure~\ref{fig:L-shape_example_deformation}.  Left: hyperelastic energy per iteration step.
   Center: maximum norm of the correction
   per iteration; Right: radius of the trust-region per iteration.  The vertical axis has logarithmic scale
   in all three images.  Note how the solver enters quadratic
   convergence after iteration 310, and how the trust-region opens up simultaneously.
   }
 \label{fig:solver_features}
\end{figure}

We start the trust-region solver at the initial iterate given in~\eqref{eq:l-shape-initial-iterate} with an initial trust-region radius of~0.1.%
\footnote{Note that this radius bounds both corrections to $m$ and to $\overline{R}$, so it cannot be
assigned a unit.}
We terminate the iteration as soon as the maximum norm of the correction drops below $3\cdot 10^{-6}$.
This criterion was achieved after 334 iterations.
Figure~\ref{fig:solver_features}, left, shows the energy $I$ per iteration (in a semi-logarithmic plot), and we observe that
the trust-region method really is monotonically energy-decreasing.  The sharp drop in the first few steps corresponds to
a decrease of the membrane energy, which dominates the initial configuration~\eqref{eq:l-shape-initial-iterate}.

Figure~\ref{fig:solver_features} also shows the correction step length and the trust-region radius
per iteration step. We note that both remain bounded in the one-digit range until the solver reaches the vicinity of the
minimizer at about iteration 310.  At this point the behavior is as predicted by Theorem~\ref{thm:tr_convergence}:
The quadratic models start to match the energy functional very well.  Correspondingly, the trust-region radius starts to
increase, and the method turns into a pure Newton method.  The expected fast local convergence can be observed
in the plot of the correction step length.
We stress that this solution is computed in a single loading step, i.e., without any
path-following mechanism.

\subsection{Torsion of a long elastic strip}
\label{sec:twisted_strip}

The purpose of the next numerical example is to show that, unlike, e.g., the approach in
\cite{gruttmann_wagner_meyer_wriggers:1993}, our discretization can easily handle
large rotations.  For this we simulate torsion of a long elastic strip, which we clamp
at one short end.  Using prescribed displacements, the other short edge is then rotated around
the center line of the strip, to a final position of three full revolutions.

\begin{table}
\begin{center}
 \begin{tabular}{cccccc}
  $h\,[\text{mm}]$ 
  & $\mu\,[\text{N}/\text{m}^2]$ & $\lambda\,[\text{N}/\text{m}^2]$ 
  & $\mu_c\,[\text{N}/\text{m}^2]$ 
  & $L_c\,[\text{mm}]$   
  & $q\,[1]$     
  \\
  \hline
  2 & $5.6452\cdot 10^9$ & $2.1796 \cdot 10^9$ & 0 & $2\cdot 10^{-3}$ & 2
 \end{tabular}
\end{center}
\caption{Material parameters for the twisted strip}
\label{tab:material_parameters_torsion}
\end{table}

Let $\omega = (0,100)\,\text{mm} \times (-5,5)\,\text{mm}$ be the parameter domain, and $\gamma_0$ and
$\gamma_1$ be the two short ends.  We clamp the shell on $\gamma_0$ by requiring
\begin{equation*}
 m(x,y) = (x,y,0),
 \qquad
 \overline{R}_3 = (0,0,1)^T
 \qquad
 \text{on $\gamma_0$},
\end{equation*}
and we prescribe a parameter dependent displacement
\begin{equation*}
 m_t(x,y) =
 \begin{pmatrix}
  1 & 0 & 0 \\
  0 & \cos 2\pi t & -\sin 2\pi t \\
  0 & \sin 2\pi t &  \cos 2\pi t
 \end{pmatrix}
 \begin{pmatrix}
  x \\ y \\ 0
 \end{pmatrix}
 \qquad
 (\overline{R}_t)_3 =
 \begin{pmatrix}
  0 \\ -\sin 2\pi t \\ \cos 2\pi t
 \end{pmatrix}
 \qquad
 \text{on $\gamma_1$}.
\end{equation*}
For each increase of $t$ by 1 this models one full revolution of $\gamma_1$ around the shell central axis.
Homogeneous Neumann boundary conditions are applied to the remaining boundary degrees of freedom.
The material parameters are
given in Table~\ref{tab:material_parameters_torsion}.  We discretize the domain with $10 \times 1$
quadrilateral elements, and use second-order (9-node) geodesic finite elements to discretize the problem.

\begin{figure}
\begin{center}
 \begin{overpic}[width=0.32\textwidth]{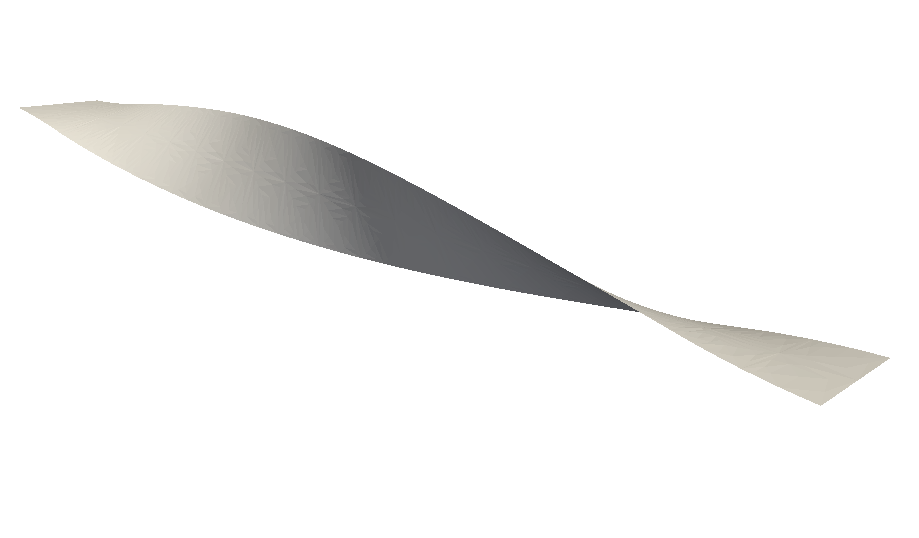}
  \put(10,15){$t = 0.5$}
 \end{overpic}
 \begin{overpic}[width=0.32\textwidth]{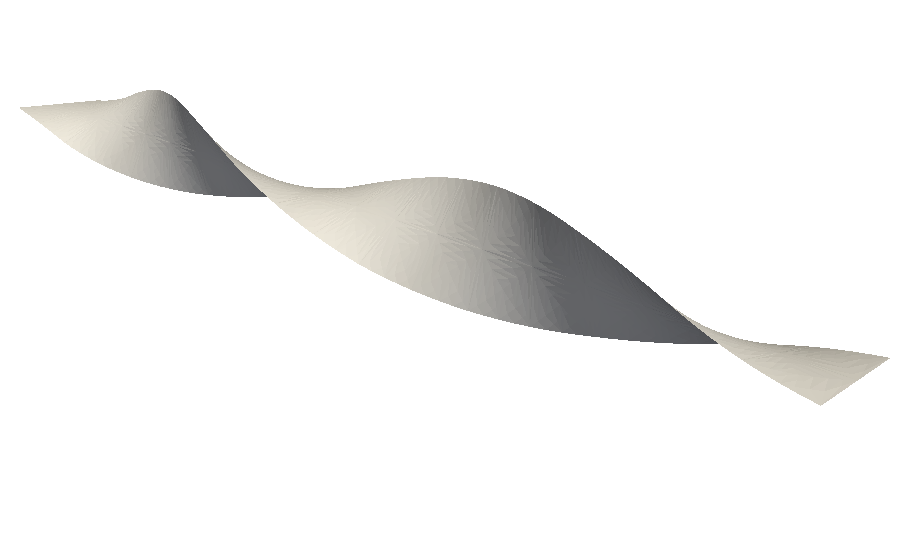}
  \put(10,15){$t = 1.0$}
 \end{overpic}
 \begin{overpic}[width=0.32\textwidth]{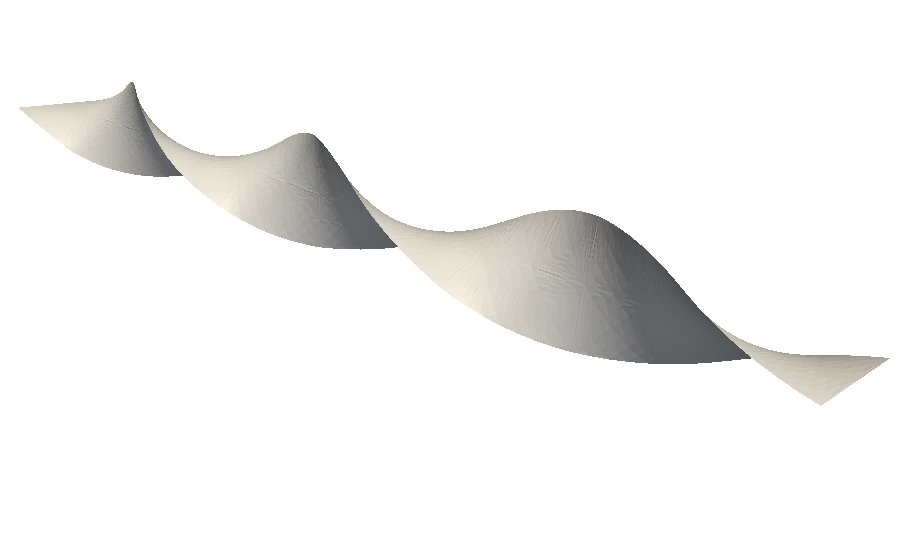}
  \put(10,15){$t = 1.5$}
 \end{overpic}

 \begin{overpic}[width=0.32\textwidth]{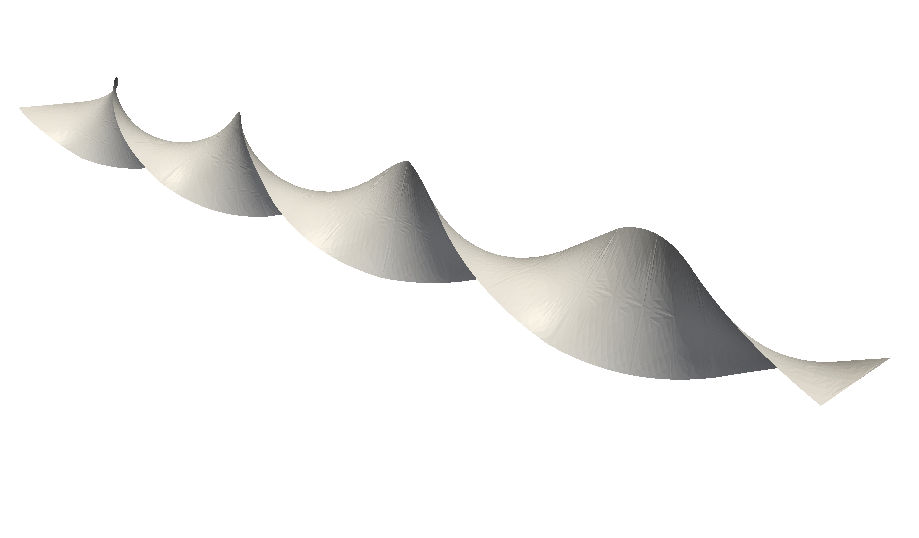}
  \put(10,15){$t = 2.0$}
 \end{overpic}
 \begin{overpic}[width=0.32\textwidth]{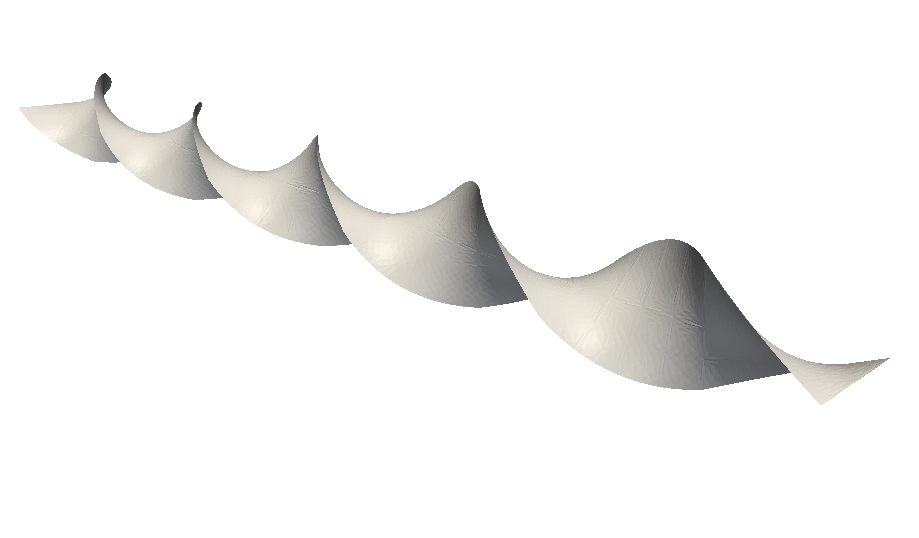}
  \put(10,15){$t = 2.5$}
 \end{overpic}
 \begin{overpic}[width=0.32\textwidth]{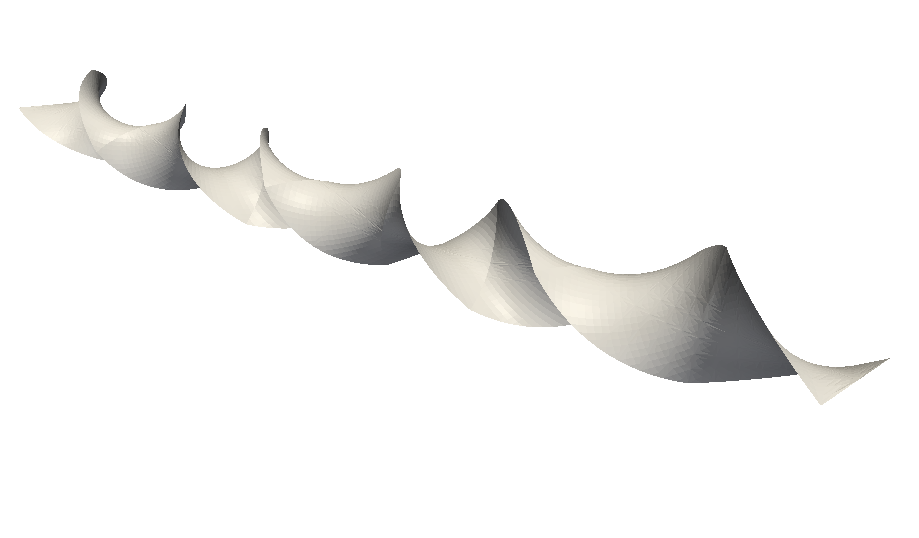}
  \put(10,15){$t = 3.0$}
 \end{overpic}
\end{center}
\caption{Twisted rectangular strip at different parameter values $t$, with $t$ equal to the number
         of revolutions.}
\label{fig:torsion_strip_solution}
\end{figure}

The result is pictured in Figure~\ref{fig:torsion_strip_solution} for several values of $t$.
Having little bending stiffness, the configuration stays symmetric throughout the parameter range.
Indeed, by increasing the length scale parameter $L_c$ one can produce materials
that are stiffer in bending.  Strips of such material buckle sideways even at only two revolutions.

In order to arrive at configurations with more than one full twist,
several intermediate loading steps have to be taken.
This is not because the Riemannian trust-region solver would not converge for $t\ge 1$.  Rather, it would
converge, but to a minimizer in the wrong homotopy group (i.e., the minimizing configuration would
never show more than a single twist).
We note also that the finite-strain membrane energy~\eqref{eq:finite_strain_membrane_energy} is essential
for this example.  Indeed, there appears to be no stable local minimizer of the small-strain
energy~\eqref{eq:small_strain_energy} that corresponds to a two-fold rotated strip.  When the energy-minimizing Riemannian
trust-region algorithm is used to minimize the small-strain energy starting from the two-revolutions
configuration, the algorithm converges to the completely planar configuration.

\subsection{Wrinkling of a sheared rectangular plastic sheet}
\label{sec:wrinkling}

In our last numerical example we demonstrate that our shell model does indeed display
microstructure.  We do this by simulating the wrinkling of a thin rectangular plastic
sheet under shearing.  Such wrinkling has been studied experimentally by
\citet{wong_pellegrino:2006a}.  Numerical simulations of their experiments can be found in~\cite{wong_pellegrino:2006c}
using the commercial FE software Abaqus, and in~\cite{taylor_bertoldi_steigmann:2014} using
a Koiter model with a finite difference discretization.
We obtain a good match between their experimental and our numerical results.

\begin{figure}
 \begin{center}
  \includegraphics[height=0.12\textheight]{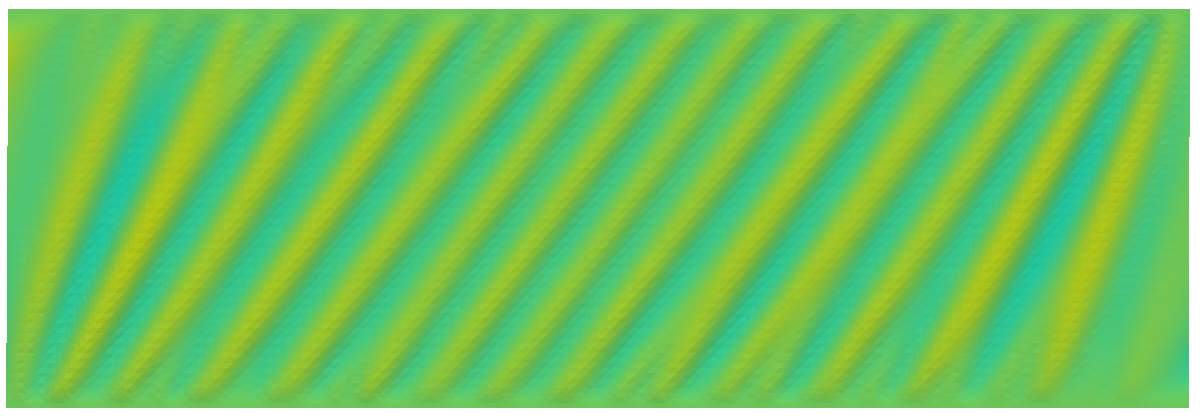}
  \hspace{0.016\textwidth}
  \begin{overpic}[height=0.12\textheight]{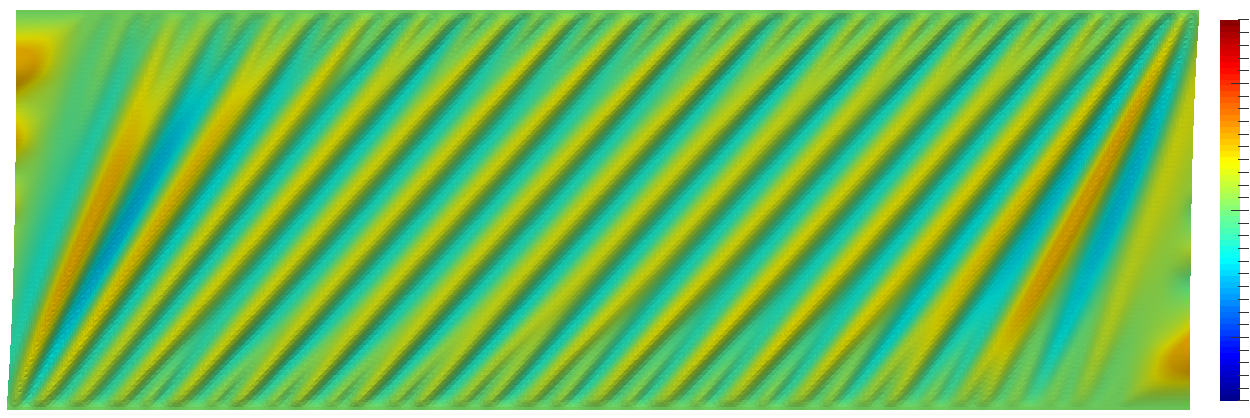}
  \put(98,35){\small $z/380$\,mm}
  \put(104,30){\small $0.003$}
  \put(101,1){\small $-0.003$}
  \put(104,15){\small $0$}
  \end{overpic}
 \end{center}
 \caption{Simulation results of the shearing tests.  The color visualizes the elevation of the wrinkles, and the color scale
   has been chosen to match the one used in \cite{taylor_bertoldi_steigmann:2014}
 }
 \label{fig:wrinkles_simulation}
\end{figure}

\begin{figure}
 \begin{center}
  \includegraphics[width=0.45\textwidth]{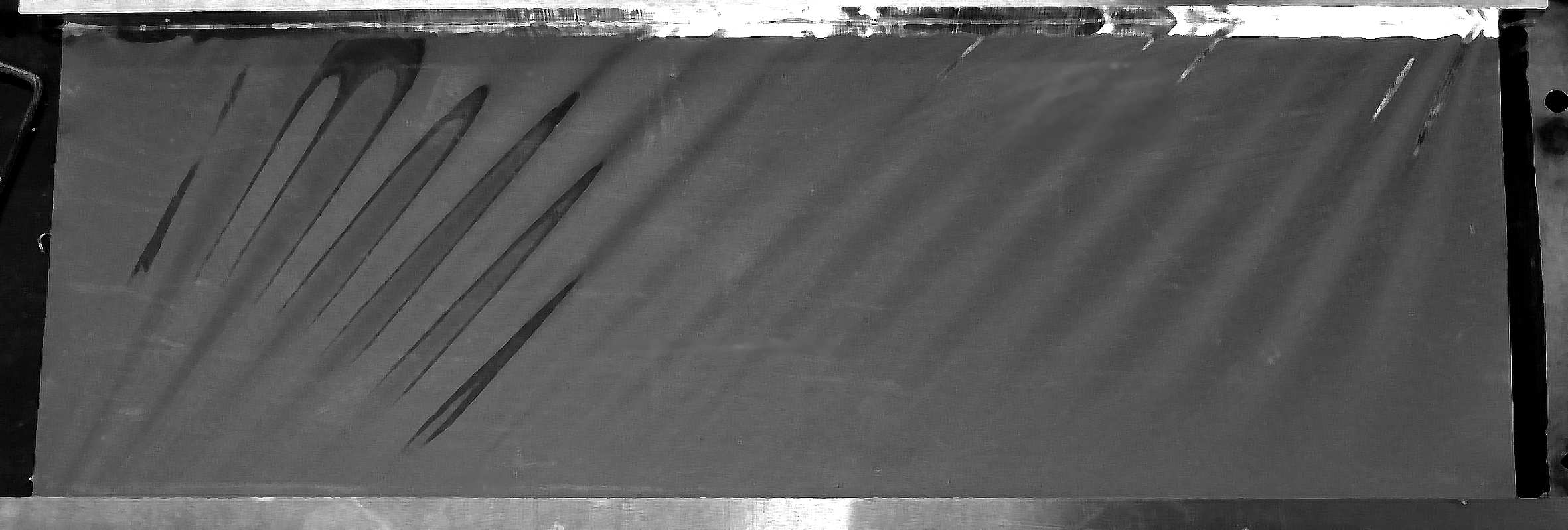}
  \hspace{0.02\textwidth}
  \includegraphics[width=0.45\textwidth]{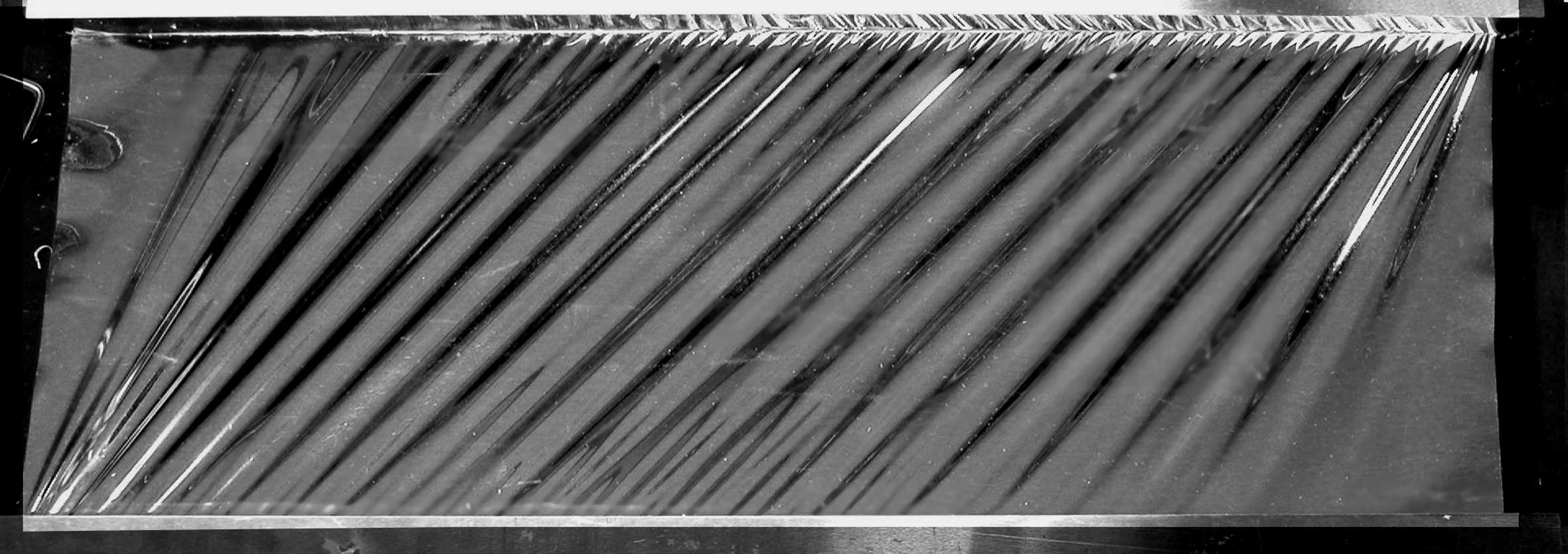}
 \end{center}
 \caption{Experimental results of the shearing tests.  Images taken from \citet{wong_pellegrino:2006a}.
 }
 \label{fig:wrinkles_experiment}
\end{figure}

The experiment consists of a rectangular plastic sheet of dimension $380\,\text{mm} \times 128\,\text{mm}$.
The sheet is clamped on the long horizontal
edges, and free on the short vertical ones.  More mathematically, we prescribe
Dirichlet boundary conditions $m(x,y) = (x,y,0)$, $\overline{R}_3(x) = (0,0,1)^T$ on the lower horizontal edge.
On the vertical sides of the domain
we prescribe zero forces and moments.  On the top horizontal side we apply
a small horizontal shearing $\delta_h$ and a vertical prestress $\delta_v$ by prescribing the Dirichlet boundary condition
$m(x,y) = (x + \delta_h, y+\delta_v, 0)$, $\overline{R}_3(x,y) = (0,0,1)^T$.

Following \citeauthor{wong_pellegrino:2006a}, we set the Lam\'e constants to
$\mu = 5.6452 \cdot 10^9\,\mathrm{N}/\mathrm{m}^2$ and $\lambda = 2.1796 \cdot 10^9\,\mathrm{N}/\mathrm{m}^2$, which corresponds
to the values $E = 3.5\,\text{GPa}$, $\nu = 0.31$ given in \cite{wong_pellegrino:2006a}.
The shell thickness is $h = 25\,\mu$m.
Additionally, we set the Cosserat couple modulus $\mu_c = 0$, the curvature exponent
$q = 2$, and the internal length scale $L_c = 0.025\,\mu\mathrm{m}$.
In \cite{wong_pellegrino:2006a}, \citeauthor{wong_pellegrino:2006a} state that they vertically prestress their sheets
slightly, but no numbers are given.  For their own numerical simulations described in \cite{wong_pellegrino:2006c}, they
use a value of $\delta_v = 0.5$\,mm.  In our own numerical experiments we found that $\delta_v = 0.5$\,mm leads to wrinkles
that are too vertical, in particular if there is not much shearing.  Low values of $\delta_v$ on the other hand do not
produce enough wrinkles.  Best results were obtained using values between $0.2$\,mm and $0.4$\,mm.

We numerically reproduce two of the four shearing experiments described in~\cite{wong_pellegrino:2006a}.  The first has a
shearing value of $\delta_h = 0.5$\,mm.  For this we discretize the domain by a structured grid with $120 \times 40 = 4\,800$
quadrilateral elements, and second-order geodesic finite elements. We set the vertical prestress to $\delta_v = 0.2$\,mm,
and start the trust-region solver from the node-wise interpolant of the function
\begin{equation*}
m(x,y) = \big(x + \delta_h y / 128\,\mathrm{mm},\; y,\; 2\,\mathrm{mm} \cos(10 x)\big),
\qquad
\overline{R}(x,y) = \operatorname{Id},
\end{equation*}
together with the Dirichlet boundary values on the top horizontal side.
The cosine waves were added to break the initial symmetry.  No attempt was made to influence the simulation results
by deliberate adjustments of the initial value.

\begin{figure}
 \begin{center}
  \includegraphics[width=0.48\textwidth]{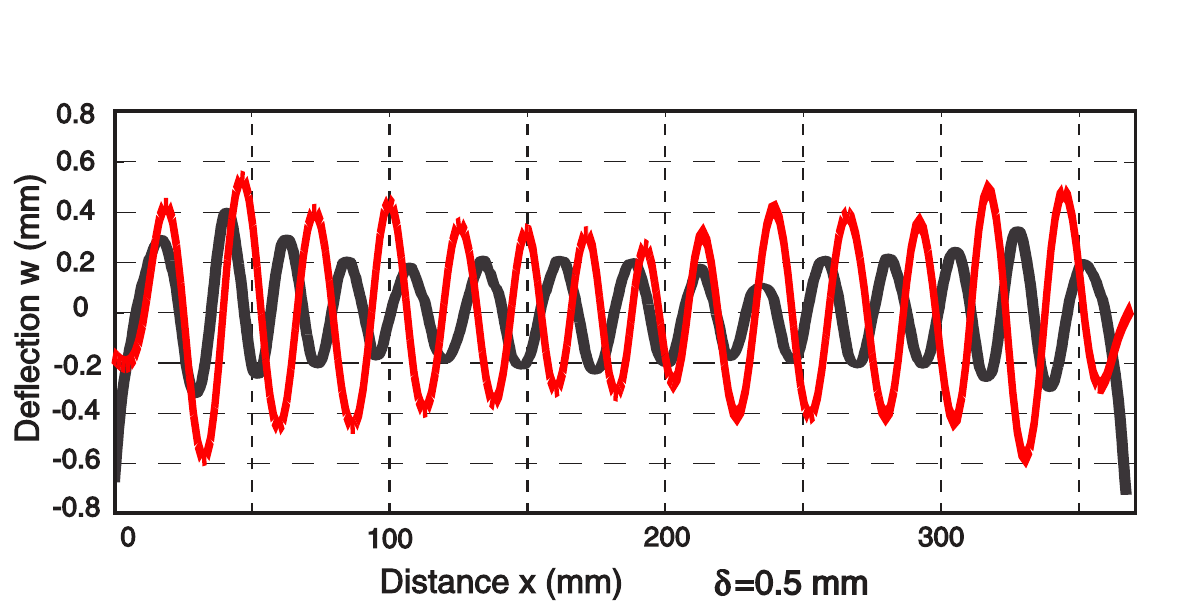}
  \hspace{0.02\textwidth}
  \includegraphics[width=0.48\textwidth]{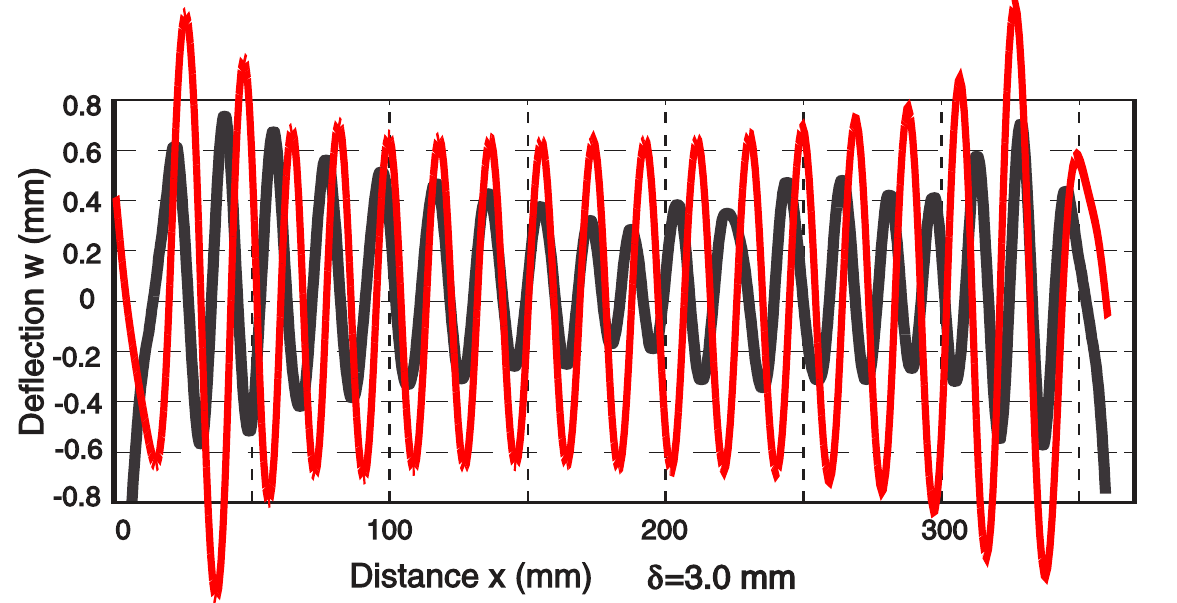}
 \end{center}
 \caption{Wrinkle amplitudes at the plane $y= 64$\,mm.  Black lines: experimental results from \citet{wong_pellegrino:2006a}.
    Red lines: our simulation results.  Observe that the number of wrinkles is almost identical, but the amplitudes
    predicted by our simulation are generally too large.
 }
 \label{fig:wrinkles_transverse_cut}
\end{figure}

Plots of the wrinkle elevation are shown on the left of Figure~\ref{fig:wrinkles_simulation}.
The results of the corresponding experiment of \citeauthor{wong_pellegrino:2006a} can be seen in Figure~\ref{fig:wrinkles_experiment},
also on the left.
We obtain a very good quantitative match with our simulation.  In particular, we obtain almost the same number
of wrinkles (Figure~\ref{fig:wrinkles_transverse_cut}).  Moreover, observe how the simulation faithfully reproduces
a lot of the fine structure, such as the secondary wrinkles near the horizontal sides, and the wrinkles
near the vertical sides.

On the other hand, the amplitudes predicted by our simulation are slightly larger than the ones observed in
the experiments.  Also, the wrinkles are inclined at a slightly steeper angle than the experimental ones..
This suggests that the prestress values $\delta_v$ is still too large.  However, as mentioned above,
a lower value of $\delta_v$ leads to a lower number of wrinkles.

The second experiment uses a larger shear value of $\delta_h = 3$\,mm.   With the other parameters as above we obtain
a result that is qualitatively correct, but the number of wrinkles is less than what \citeauthor{wong_pellegrino:2006a}
observed in their experiments.  A better match is obtained by increasing the vertical prestress to $\delta_v = 0.4$\,mm
and using a fine grid with $240 \times 80 = 19\,200$ elements.  This simulation is what is plotted on the right of
Figures~\ref{fig:wrinkles_simulation}, \ref{fig:wrinkles_experiment}, and \ref{fig:wrinkles_transverse_cut}.
Now we observe a very good quantitative agreement also for this more extreme case, with the same restrictions
as for the low-shear case.  Since we have not observed artificial stiffness introduced by our discretization,
we suspect that using the finer grid makes the trust-region algorithm end up in a different local minimizers of the energy.

\section{Appendix: Implementation of geodesic finite elements for \texorpdfstring{$\text{SO}(3)$}{SO(3)}}

In this appendix we explain how the geodesic interpolation (Definition~\ref{def:geodesic_interpolation})
that forms the basis of the geodesic finite element method
can be implemented in practice.  Since the definition of the interpolation function
\begin{align}
 \nonumber
 \Upsilon & \; : \; \SOdrei^m \times T_\text{ref} \to \SOdrei \\
 \label{eq:geodesic_interpolation_revisited}
 \Upsilon(R_1,\dots,R_m;\xi) & = \argmin_{Q \in \SOdrei} \sum_{i=1}^m \lambda_i(\xi) \dist(R_i,Q)^2
\end{align}
uses a minimization formulation, its use in practice warrants a few explanations.

The Cosserat shell energies of Chapter~\ref{sec:continuous_model} are both
first-order energies.  Hence, to evaluate them for a given geodesic finite element function $\overline{R}_h$ we need to compute
function values $\overline{R}_h(x) \in \SOdrei$ and first derivatives $\nabla \overline{R}_h(x) : \R^2 \to T_{\overline{R}_h(x)}\SOdrei$
at given (quadrature) points $x \in \omega$.
Using the integral transformation formula this can be reduced to computing values and
first derivatives of the interpolation function $\Upsilon$ on the reference element $T_\text{ref}$.

Finding minimizers of the energy by a Riemannian trust-region method additionally requires the gradient $\nabla \bar{I}$
and the Hessian $\hessian \bar{I}$ of the algebraic Cosserat shell energy~\eqref{eq:algebraic_shell_energy}.
By the chain rule, expressions for these include derivatives of $\Upsilon$ and $\nabla \Upsilon$ with respect
to the coefficients $R_1,\dots,R_m$.  These can in principle be computed semi-analytically~\cite{sander:2012}.
However, we have found using an automatic differentiation system
much more convenient (see Section~\ref{sec:computing_tangent_problem}).

\subsection{Quaternion coordinates for \texorpdfstring{$\text{SO}(3)$}{SO(3)}}
\label{sec:quaternion_coordinates}

While the construction and theory of geodesic finite elements is completely coordinate-free,
an implementation necessarily needs some sort of coordinates for $\SOdrei$.
The naive approach uses the canonical embedding of $\SOdrei$ into $\R^{3 \times 3}$.
However, quaternion coordinates allow a more efficient implementation.

Let
\begin{equation*}
 \mathbb{H}_{\abs{1}} \colonequals \big\{ p \in \R^4 \; \big| \; \abs{p} = 1 \big\}
\end{equation*}
be the set of unit quaternions, i.e., the unit sphere $S^3 \subset \R^4$
equipped with quaternion multiplication
\begin{equation*}
\cdot : \mathbb{H}_{\abs{1}} \times \mathbb{H}_{\abs{1}} \to \mathbb{H}_{\abs{1}}
\qquad\qquad
 p \cdot q
 =
 \begin{pmatrix}
   p_3 q_0 - p_2 q_1 + p_1 q_2 + p_0 q_3 \\
   p_2 q_0 + p_3 q_1 - p_0 q_2 + p_1 q_3 \\
 - p_1 q_0 + p_0 q_1 + p_3 q_2 + p_2 q_3 \\
 - p_0 q_0 - p_1 q_1 - p_2 q_2 + p_3 q_3
 \end{pmatrix}
.
\end{equation*}
The unit quaternions $\mathbb{H}_{\abs{1}}$ form a smooth compact manifold embedded in $\R^4$,
and global coordinates on $\mathbb{H}_{\abs{1}}$ are naturally given by this embedding.
The tangent space at a point $p \in \mathbb{H}_{\abs{1}}$ is
\begin{equation*}
 T_p \mathbb{H}_{\abs{1}} = \{ v \in \R^4 \; | \; \langle p,v \rangle_{\R^4} = 0 \},
\end{equation*}
hence tangent vectors $v \in T_p \mathbb{H}_{\abs{1}}$ can be treated as vectors in $\R^4$.
For any $q\in \mathbb{H}_\abs{1}$, the projection $P_q : T_q\R^4 \to T_q \mathbb{H}_\abs{1}$ is given by
\begin{equation*}
 (P_q)_{ij} = \delta_{ij} - q_i q_j.
\end{equation*}

A Riemannian structure for $\mathbb{H}_{\abs{1}}$ is obtained by inheriting the metric of the surrounding space
\begin{equation*}
 \langle v,w \rangle_{T_p \mathbb{H}_{\abs{1}}}
\colonequals
 \langle v,w \rangle_{\R^4}
\qquad
\text{for all $v,w \in T_p \mathbb{H}_{\abs{1}}$}.
\end{equation*}
For a point $p \in \mathbb{H}_{\abs{1}}$ and a tangent vector $v \in T_p \mathbb{H}_{\abs{1}}$, the
exponential map $\exp_p : T_p \mathbb{H}_{\abs{1}} \to \mathbb{H}_{\abs{1}}$ is then given by
\cite[Ex.\,5.4.1]{absil_mahony_sepulchre:2008}
\begin{equation}
\label{eq:quaternion_exponential_map}
 \exp_p v = \cos \abs{v} \cdot p + \frac{\sin \abs{v}}{\abs{v}} \cdot v.
\end{equation}

The unit quaternions can be used to represent rotations, because there is a natural relationship between $\mathbb{H}_{\abs{1}}$
and $\SOdrei$.  More precisely, the map $F : \mathbb{H}_{\abs{1}} \to \R^{3 \times 3}$
\begin{equation}
\label{eq:SO3_double_cover}
 F(p)
 \colonequals
 \begin{pmatrix}
  p_0^2 - p_1^2 - p_2^2 + p_3^2  &  2(p_0 p_1 - p_2 p_3)               & 2 (p_0 p_2 + p_1 p_3) \\
  2(p_0 p_1 + p_2 p_3)           &  - p_0^2 + p_1^2 - p_2^2 + p_3^2    & 2 (p_1 p_2 - p_0 p_3) \\
  2(p_0 p_2 - p_1 p_3)           &  2(p_0 p_3 + p_1 p_2)               & - p_0^2 - p_1^2 + p_2^2 + p_3^2
 \end{pmatrix}
\end{equation}
is a Lie group homomorphism from $\mathbb{H}_{\abs{1}}$ onto $\SOdrei$.  It is two-to-one, meaning that for each
point $p \in \mathbb{H}_{\abs{1}}$ there is exactly one other point, namely $-p$, representing the same rotation
$F(p) = F(-p) \in \SOdrei$.  Using quaternion coordinates for rotations reduces the memory footprint and computing times
considerably.  For the rest of this chapter we use upper case letters $Q,R$ for elements of $\SOdrei$, and lower case
letters $p,q$ for quaternions.

\subsection{The canonical distance of \texorpdfstring{$\text{SO}(3)$}{SO(3)} in quaternion coordinates}

The metric structure of the set of unit quaternions is identical to metric structure of the unit sphere in $\R^4$.
The geodesics of $\mathbb{H}_{\abs{1}}$ are the segments of great circles.  Any two points $p,q \in \mathbb{H}_{\abs{1}}$ can
be connected by such segments; hence $\mathbb{H}_{\abs{1}}$ is geodesically complete.
If $p \ne -q$ there is a unique shortest geodesic that connects $p$ and $q$.
For all pairs of points $p = -q$ there are infinitely many minimizing
geodesics, each of length $\pi$.  Hence the injectivity radius of $\mathbb{H}_{\abs{1}}$
is $\inj (\mathbb{H}_{\abs{1}}) = \pi$.

The Riemannian distance between two points $p$ and $q$ is the length of the
shortest arc of a great circle connecting $p$ to $q$.  Let $\gamma : [0,1] \to S^3$ be such an arc.
Its length is given by
\begin{equation}
 \label{eq:distance_on_sm}
 \dist_{\mathbb{H}_{\abs{1}}}(p,q)
 =
 \int_0^1 \abs{\gamma'(t)}_{S^3}\,dt
 =
 \arccos \langle p,q \rangle_{\R^4}.
\end{equation}

We now use this to express the canonical distance on $\text{SO}(3)$ in terms of quaternion coordinates.
To avoid confusion we now always write $\dist_{\mathbb{H}_{\abs{1}}}$ or $\dist_{\SOdrei}$.
First note that $F$ defined in \eqref{eq:SO3_double_cover} is a scaling in the sense that
\begin{equation}
\label{eq:scaling_of_F}
 \norm{\nabla F \cdot v}_{\R^{3\times 3}} = 2 \norm{v}_{\R^4}
\end{equation}
for any $v \in T_q \mathbb{H}_{\abs{1}}$.
Let $R_1, R_2 \in \SOdrei$ be two rotations and let $p,q \in \mathbb{H}_{\abs{1}}$ be such that $F(p) = R_1$
and $F(q) = R_2$.  We first consider the simpler case that $\dist_{\mathbb{H}_{\abs{1}}}(p,q) < \pi/2$.
Suppose that $\gamma$ is the shortest path from $p$ to $q$.  Then, by \eqref{eq:scaling_of_F}, $F(\gamma)$ is a shortest path from $F(p)$ to $F(q)$
in $\SOdrei$, and
\begin{equation*}
 \dist_{\SOdrei}(F(p), F(q))
 =
 \int_0^1 \abs{(F \circ \gamma)'(t)}_{\SOdrei}\,dt
  =
 2 \int_0^1 \abs{\gamma'(t)}_{S^3}\,dt
 =
 2 \dist_{\mathbb{H}_{\abs{1}}}(p,q).
\end{equation*}

For the general case, we also have to take into account that $\mathbb{H}_{\abs{1}}$ is a double cover of $\text{SO}(3)$.
Let $p,q \in \mathbb{H}_{\abs{1}}$ be such that $\dist_{\mathbb{H}_{\abs{1}}}(p,q) > \pi/2$.
Then $q$ represents the same element of $\text{SO}(3)$ as $-q$, but
$\dist_{\mathbb{H}_{\abs{1}}}(p,-q) = \pi - \dist_{\mathbb{H}_{\abs{1}}}(p,q) < \pi/2$.
The distance on $\text{SO}(3)$ for arbitrary $p,q$ given in terms
of the distance on $\mathbb{H}_{\abs{1}}$ is therefore
\begin{equation}
\label{eq:dist_SO_quat_coordinates}
 \dist_{\SOdrei}(F(p),F(q))
 =
 \begin{cases}
  2\dist_{\mathbb{H}_{\abs{1}}}(p,q) & \text{if $\dist_{\mathbb{H}_{\abs{1}}}(p,q) \le \pi/2$}, \\
  2\pi - 2\dist_{\mathbb{H}_{\abs{1}}}(p,q) & \text{otherwise}.
 \end{cases}
\end{equation}
Note that this metric is continuous, but not differentiable at points $p,q$ with $\dist_{\mathbb{H}_{\abs{1}}}(p,q) = \pi/2$.
This comes as no surprise as this is precisely the case when $F(q)$ is in the cut locus of $F(p)$.

\bigskip

For an algorithmic evaluation of the interpolation formula~\eqref{eq:geodesic_interpolation_revisited} we will need first and second derivatives
of $\dist_{\SOdrei}(R,\cdot)^2$
with respect to its second argument, for fixed arbitrary $R \in \SOdrei$.  We use~\eqref{eq:dist_SO_quat_coordinates}, and Lemmas~\ref{lem:embedded_gradient} and~\ref{lem:embedded_hessian} on the derivatives
of scalar-valued functions on embedded manifolds.
For these, we need an extension of $\dist_{\mathbb{H}_{\abs{1}}}$ to a neighborhood of $\mathbb{H}_{\abs{1}}$ in $\R^4$.
We choose
\begin{equation*}
 \widetilde{\dist}_{\mathbb{H}_\abs{1}}(p,q)^2
   \colonequals \dist_{\mathbb{H}_\abs{1}}\Big(p, \frac{q}{\abs{q}}\Big)^2
   = \arccos^2 \Big\langle p, \frac{q}{\abs{q}}\Big\rangle.
\end{equation*}
This is well-defined and smooth on a neighborhood of $\mathbb{H}_\abs{1}$ in $\R^4$.
For ease of notation we define $\alpha : [-1,1] \to \R$, $\alpha(x) \colonequals \arccos^2(x)$.

We now compute the first derivative of $\dist_{\SOdrei}(R,\cdot)^2$
\begin{equation}
\label{eq:derivative_of_dist_no_coefficients}
 \parder{}{q} \dist_{\SOdrei}(R,F(q))^2
  \in T_q \mathbb{H}_\abs{1},
\end{equation}
for arbitrary but fixed $R \in \SOdrei$.  Note that
\begin{equation*}
 \parder{}{q} \Big\langle p, \frac{q}{\abs{q}}\Big\rangle = P_q p
 \qquad
 \text{if $\abs{q} = 1$.  }
\end{equation*}
With \eqref{eq:projected_gradient}, \eqref{eq:distance_on_sm}, and $\abs{q} = 1$ we get for the coefficients $i=1,\dots,4$ of \eqref{eq:derivative_of_dist_no_coefficients}
\begin{align*}
 \parder{}{q_i} \dist_{\SOdrei}(R,F(q))^2
 & = \Big(\parder{}{q} \dist_{\SOdrei}(R,F(q))^2 \Big)_i \\
 & =
 \begin{cases}
   4\alpha'(x)\Big|_{x=\langle p,q \rangle} (P_qp)_i & \text{if $\dist_{\mathbb{H}_{\abs{1}}} (p,q) \le \pi/2$} \\
   -4\alpha'(x)\Big|_{x=\langle p,q \rangle} (P_qp)_i & \text{else},
 \end{cases}
\end{align*}
where $p$ is any one of the two points on $\mathbb{H}_{\abs{1}}$ with $F(p) = R$.
Note that since $\dist_{\mathbb{H}_{\abs{1}}}(p,q) \le \pi/2$ if and only if $\langle p,q\rangle \ge 0$ this is equivalent to
\begin{equation}
\label{eq:first_derivative_of_dist}
 \parder{}{q_i} \dist(p,q)^2
=
  \operatorname{sgn}\big[\langle p,q\rangle\big] 4\alpha'(x)\Big|_{x=\abs{\langle p,q \rangle}} (P_qp)_i.
\end{equation}
The derivative of $\alpha(x) = \arccos^2(x)$ can be given in closed form
\begin{equation*}
\alpha'(x) = - \frac{2 \arccos(x)}{\sqrt{1-x^2}}.
\end{equation*}
However, this expression gets numerically unstable around $x=1$.  There, the series expansion
\begin{equation*}
\alpha'(x)
 =
-2+ \frac{2 (x-1)}{3}
        + O((x-1)^2)
\end{equation*}
has to be used instead.

For the second derivative of $\dist_{\SOdrei}(R,\cdot)^2$ we note that
\begin{equation*}
 \parder{}{p_j} (P_q p)_i = (P_q)_{ij} = \delta_{ij} - q_i q_j
\qquad \text{and} \qquad
 \parder{}{q_j} (P_q p)_i = - \delta_{ij} \langle p,q \rangle -  q_i p_j
\end{equation*}
for any $p \in \R^4$.
Using Lemma~\ref{lem:embedded_hessian} we obtain
\begin{equation}
\label{eq:second_derivative_of_dist}
 \Big[\parder{^2}{q^2} \dist_{\SOdrei}(R,F(q))^2 \Big]_{ij}
=
  4\alpha'' (P_qp)_i (P_qp)_j
 -4\operatorname{sgn}(\langle p,q \rangle)\alpha' (P_q)_{ij} \langle p,q \rangle,
\end{equation}
where again $p \in \mathbb{H}_{\abs{1}}$ is such that $F(p) = R$.

The second derivative of $\alpha(x)$ is
\begin{equation*}
\alpha''(x) = (\arccos^2(x))''  = \frac{2}{1-x^2} - \frac{2x \arccos(x)} {(1-x^2)^{3/2}}.
\end{equation*}
Again, near $x=1$ this gets unstable and has to be replaced by its series expansion
\begin{equation*}
\alpha''(x)
  =
\frac{2}{3}- \frac{8}{15}(x-1) +O((x-1)^2).
\end{equation*}

\subsection{Evaluation of geodesic interpolation functions}
\label{sec:values_of_gfe_functions}

We now discuss how values and first derivatives of the interpolation function $\Upsilon$ can be computed in practice.
Unfortunately, there are no closed-form expressions for the solution of the minimization
problem~\eqref{eq:geodesic_interpolation_revisited}, and it therefore
needs to be solved numerically.  As its objective functional (written in quaternion coordinates)
\begin{equation*}
f_\xi
:
q \mapsto \sum_{i=1}^m \lambda_i(\xi) \dist_{\SOdrei}(R_i, F(q))^2
\end{equation*}
is defined on the Riemannian manifold $\mathbb{H}_{\abs{1}} \subset \R^4$ we use a Riemannian Newton method
as presented in~\cite{absil_mahony_sepulchre:2008}.%
\footnote{In~\cite{sander:2012} it was proposed to use a Riemannian trust-region method instead
of the simpler Newton method.  Such a choice guarantees convergence of the solver.  However, in practice we never observed
convergence issues even for the simpler Newton method.}
Under the assumptions of Theorems~\ref{thm:unique_minimizer_first_order} (for $p=1$) and~\ref{thm:well_posedness_vague}
(for $p>1$),
$f_\xi$ is $C^\infty$ (\cite[Lem.\,2.4]{sander:2012}), and
strictly convex on an open geodesic ball containing the $R_i$ (\cite[Thm.\,1.2]{karcher:1977}
and~\cite[Lem.\,3.11]{sander:2013}, respectively).

One step of the Riemannian Newton method on $\mathbb{H}_{\abs{1}}$ takes the following form.
With $k$ the iteration number let $q_k \in \mathbb{H}_{\abs{1}}$ be the current iterate.
We use the exponential map $\exp_{q_k} : T_{q_k} \mathbb{H}_{\abs{1}} \to \mathbb{H}_{\abs{1}}$ (see~\eqref{eq:quaternion_exponential_map}) to define lifted functionals%
\begin{equation*}
 \hat{f}_k : T_{q_k} \mathbb{H}_{\abs{1}} \to \R
 \qquad
 \hat{f}_k (s) \colonequals f_\xi (\exp_{q_k} s).
\end{equation*}
The Newton update at step $k$ is then
\begin{equation}
\label{eq:local_newton_update}
q_{k+1} = \exp_{q_k} s_k
\qquad \text{with} \qquad
s_k = - \hessian \hat{f}_k(0)^{-1} \nabla \hat{f}_k(0).
\end{equation}
Using $\nabla \exp 0 = \identity$ we see that the gradient of $\hat{f}_k$ at $0 \in T_{q_k} \mathbb{H}_{\abs{1}}$ is
\begin{equation*}
 \nabla \hat{f}_k(0)
=
\sum_{i=1}^m \lambda_i(\xi) \frac{\partial}{\partial q} \dist_{\SOdrei}(R_i, F(q))^2,
\end{equation*}
and that the Hessian is
\begin{equation*}
 \hessian \hat{f}_k(0)
 =
  \sum_{i=1}^m \lambda_i(\xi) \frac{\partial^2}{\partial q^2} \dist_{\SOdrei}(R_i, F(q))^2.
\end{equation*}
The two derivatives of the distance function have been given in~\eqref{eq:first_derivative_of_dist} and~\eqref{eq:second_derivative_of_dist}.
The matrix $\hessian \hat{f}_k(0)$ is $4 \times 4$,
and has a one-dimensional kernel, which is the normal space of $S^3$ in $\R^4$ at
$q_k$.  We use a rank-aware direct solver for the Newton update systems~\eqref{eq:local_newton_update}.
The Newton solver typically needs only a handful of iterations to converge up to machine precision.

\bigskip

In the proof of
Lemma~\ref{thm:geodesic_interpolation_is_differentiable} the implicit function theorem
was used to show under what circumstances the derivative $\partial \Upsilon / \partial \xi$ exists.
Here we use it again for the actual computation.  For ease of notation we introduce $\widetilde{\Upsilon} \colonequals F^{-1}(\Upsilon)$,
which gives interpolation points expressed as quaternions.
By \cite[Lem.\,2.4]{sander:2012} the functional $f_\xi$ is smooth.
Hence, its minimizer can be characterized by
\begin{equation}
\label{eq:average_distance_gradient_is_zero}
 \Phi(R_1,\dots,R_m; \xi, \widetilde\Upsilon(R_1,\dots,R_m;\xi)) = 0,
\end{equation}
where
\begin{align}
\nonumber
 \Phi & \; : \; \SOdrei^m \times T_\text{ref} \times \mathbb{H}_{\abs{1}} \to T\mathbb{H}_{\abs{1}} \\
\label{eq:gradient_of_coordinate_functional}
 \Phi(R_1,\dots,R_m,\xi,q)
  & =
 \parder{}{q} f_\xi(q) = \sum_{i=1}^m \lambda_i(\xi) \parder{}{q} \dist_{\SOdrei} (R_i,F(q))^2 .
\end{align}
Taking the total derivative of \eqref{eq:average_distance_gradient_is_zero}
with respect to $\xi$ we get
\begin{multline*}
\totalder{}{\xi} \Phi(R_1,\dots,R_m,\xi,\widetilde\Upsilon(R_1,\dots,R_m; \xi)) \\
   =
   \parder{\Phi(R_1,\dots,R_m,\xi,q)}{\xi} + \parder{\Phi(R_1,\dots,R_m,\xi,q)}{q} \cdot \parder{\widetilde\Upsilon(R_1,\dots,R_m;\xi)}{\xi} = 0.
\end{multline*}
By \cite[Lem.\,3.11]{sander:2013} the matrix
\begin{equation*}
\frac{\partial \Phi}{\partial q}
  =
\hessian f_\xi
 \in
 \R^{4 \times 4}
\end{equation*}
is invertible on the three-dimensional subspace $T_{\widetilde\Upsilon(R_1,\dots,R_m;\xi)} \mathbb{H}_{\abs{1}} \subset \R^4$,
and hence $\partial \widetilde\Upsilon(R_1,\dots,R_m;\xi) / \partial \xi$
can be computed as the solution of the linear system of equations
\begin{equation}
 \label{eq:linear_system_for_gradient}
 \parder{\Phi(R_1,\dots,R_m,\xi,q)}{q} \cdot \parder{\widetilde\Upsilon(R_1,\dots,R_m;\xi)}{\xi} = - \parder{\Phi(R_1,\dots,R_m,\xi,q)}{\xi}.
\end{equation}
Using the definition \eqref{eq:gradient_of_coordinate_functional} we see that
in coordinates
$\partial \Phi / \partial \xi$ is a $4 \times 2$-matrix, where the $i$-th column is
\begin{equation*}
 \Big(\parder{\Phi}{\xi}^T\Big)_i = \frac{\partial}{\partial q} \dist_{\SOdrei}(R_i,F(q))^2.
\end{equation*}
Hence evaluating the derivative of a geodesic finite element function amounts to
an evaluation of its value (to know where to evaluate the derivatives
of $\Phi$) and the solution of the
symmetric linear system~\eqref{eq:linear_system_for_gradient}.

\section*{Acknowledgements}

The authors would like to thank Kshitij Kulshreshtha for his help with the ADOL-C automatic differentiation system,
and Ingo Münch for the interesting discussions on the discretization of finite strain Cosserat problems.

\bibliography{SanderNeffBirsan_nonlinear_Cosserat_shells}

\end{document}